\documentclass[a4paper, final, 12pt, twoside]{article}
\usepackage{amsmath, amssymb, latexsym, amscd, amsthm,amsfonts,amstext}
\usepackage[mathscr]{eucal}
\usepackage{graphicx}
\usepackage{subfig}
\usepackage{float}
\usepackage{color}
\usepackage{hyperref}
\usepackage[utf8]{inputenc}
\usepackage[english]{babel}
\usepackage{fancyhdr}

\numberwithin{equation}{section}
\input{tcilatex}

\begin{document}

\title{On the travel time tomography problem in 3D }
\author{Michael V. Klibanov\thanks{%
Department of Mathematics and Statistics, University of North Carolina
Charlotte, Charlotte, NC, 28223, mklibanv@uncc.edu.} }
\date{}
\maketitle

\begin{abstract}
Numerical issues for the 3D travel time tomography problem with
non-overdetemined data are considered. Truncated Fourier series with respect
to a special orthonormal basis of functions depending on the source position
is used. In addition, truncated trigonometric Fourier series with respect to
two out of three spatial variables is used. First, the Lipschitz stability
estimate is obtained. Next, a globally convergent numerical method is
constructed using a Carleman estimate for an integral operator.
\end{abstract}

\textbf{Key Words}. global convergence, semi-finite dimensional mathematical
model, Carleman Weight Function, weighted Tikhonov-like functional

\textbf{AMS subject classification}. 35R25, 35R30

\section{Introduction}

\label{sec:1}

We call a numerical method for a nonlinear inverse problem \emph{globally}
convergent if there exists a theorem claiming that this method delivers
points in a sufficiently small neighborhood of the correct solution without
any advanced knowledge of this neighborhood. In other words, a good first
guess is not needed. We construct here a globally convergent numerical
method for the Travel Time Tomography Problem (TTTP) in the 3D isotropic
case with formally determined incomplete data. The TTTP is also called
sometimes Inverse Kinematic Problem of Seismic, see chapter 3 in \cite{R}.
Previous publications about the 3D TTTP work only with the over determined
data, in which case the number of free variables in the data $m=4$ exceeds
the number $n=3$ of free variables in the unknown coefficient, see, e.g. 
\cite{MR,R,Stef} for the theoretical results and \cite{Z} for a numerical
method. On the other hand, our data are formally determined with $m=n=3$. In
the 2D case, the TTTP\ is always formally determined with $m=n=2,$ see \cite%
{BG,Mukh,PU} for the theory and \cite{Sch} for an algorithm.

The TTTP was first considered by Herglotz \cite{Herg} and Wiechert and
Zoeppritz \cite{W} in 1905 and 1907 respectively in the 1D case, due to
important applications in Geophysics, see chapter 3 of \cite{R} for some
details. Recently it was discovered that, in addition to Geophysics, the
TTTP also arises in the inverse problem of the recovery of the spatially
distributed dielectric constant coefficient from the scattering data without
the phase information in both Helmholtz equation \cite{KR} and Maxwell
Equations \cite{RMaxwell}.

The first globally convergent numerical method for the 3D TTTP with formally
determined data was constructed in the recent work of the author \cite{KTTTP}%
. In addition, the data in \cite{KTTTP} are incomplete, so as in the current
paper, see Figure 1 in \cite{KTTTP}. The globally convergent numerical
method of this paper differs from the one of \cite{KTTTP} in one important
aspect. In \cite{KTTTP} finite differences with respect to two $x,y$ out of
three $x,y,z$ spatial variables were considered.\ On the other hand, we
consider here truncated trigonometric Fourier series with respect to $x,y$.
The replacement of those finite differences with the truncated trigonometric
Fourier series significantly simplifies and shortens the presentation.

In numerical methods for inverse problems, it is quite acceptable to use 
\emph{approximate mathematical models}, see, e.g. \cite{BK,GN,Kab0,Kab1,Kab2}%
. We work below only within the framework of such a model. Our approximate
mathematical model consists of two assumptions. First, we assume that a
certain function associated with the solution of the governing eikonal
equation can be represented via a truncated Fourier series with respect to a
special orthonormal basis in $L_{2}\left( 0,\pi \right) $. This basis was
recently proposed by the author in \cite{KJIIP}. Functions of that basis
depend only on the position of the point source. Second, we assume, as
stated above, that each component of that Fourier series can be represented
via a truncated trigonometric Fourier series with respect to $x,y.$ A
different approximate mathematical model for the TTTP is in \cite{KTTTP},
where finite differences instead of these truncated trigonometric Fourier
series are used.

Thus, we work here with a semi-finite dimensional approximation of the
original inverse problem. Assumptions of this sort are often used in
numerical methods for inverse problems.\ Furthermore, convergence when the
number of terms of truncated series tends to infinity is usually not proven,
see, e.g. \cite{GN,Kab0,Kab1,Kab2}. Indeed, it is well known that such
convergence results are very challenging to prove due to the ill-posed
nature of inverse problems.

A conventional Tikhonov least squares cost functional for a nonlinear
inverse problem is non convex. Hence, there is no guarantee that this
functional does not have local minima and ravines, see, e.g. \cite{Scales}
for a good numerical example of local minima even for a rather simple
inverse problem. Local minima might lead to incorrect solutions since any
gradient-like method of the optimization of that functional can stop at any
local minimum. To avoid the local minima, it was proposed in \cite%
{BK1,Klib97,KT,KK} to construct globally strictly convex cost functionals
for Coefficient Inverse Problems (CIPs) and in \cite{KQ} for ill-posed
problems for quasilinear PDEs. This procedure is called \emph{convexification%
}. Next, it was established in \cite{Bak} that the minimizer of such a
functional exists, is unique and minimizers converge to the exact solution
of the original problem as long as the noise in the data tends to zero.
Finally, it was also established in \cite{Bak} that the gradient projection
method converges to the exact solution if starting from an arbitrary point
of a certain ball in the Hilbert space of an arbitrary radius $R>0.$ Since
smallness conditions are not imposed on $R$, then this is the \emph{global}
convergence by our above definition. After the theory was cleared in \cite%
{Bak}, a number of works were published in which numerical studies of the
convexification were conducted \cite{KlibKol1,convnew,KlibKol2,KEIT,Ktime}.
We also refer here to the publication \cite{Baud} for a different version of
the convexification for a hyperbolic coefficient inverse problem with a non
vanishing initial condition.

The main ingredient of the convexification is the presence of the Carleman
Weight Function(CWF) in the weighted Tikhonov-like functional. CWF is the
function which is involved as the weight function in the Carleman estimate
for a corresponding PDE operator. Unlike this, in the current paper, so as
in \cite{KTTTP}, we apply the CWF for a Volterra linear integral operator.
As a result, we construct here a globally strictly convex weighted
Tikhonov-like functional for the TTTP. Note that prior \cite{KTTTP} CWFs
were used only for PDE operators. Also, \cite{KTTTP} is the first work in
which a Carleman estimate is applied to the TTTP.

In section 2 we pose the inverse problem. In section 3 we present some
preliminary considerations. In particular, we provide more insight than
before in our orthonormal basis of \cite{KJIIP}. In section 4 we formulate
our semi-finite dimensional approximation for a certain function we work
with and also describe our approximate mathematical model. In section 5 we
prove Lipschitz stability result for our approximate mathematical model. In
section 6 we describe our numerical method and formulate corresponding
theorems. We prove these theorems in section 7.

\section{Statement of the Problem}

\label{sec:2}

Below $\mathbf{u}=\left( x,y,z\right) $ denotes points in $\mathbb{R}^{3}.$
Let $B_{1},B_{2}>0$ be two numbers, $B_{2}>B_{1}$. Define the domain $%
G\subset \mathbb{R}^{3}$ as%
\begin{equation}
G=\left\{ \mathbf{u}=\left( x,y,z\right) :0<x,y<\pi ,z\in \left(
B_{1},B_{2}\right) \right\} .  \label{2.1}
\end{equation}%
The lower and upper boundaries of $G$ are denoted as $\Gamma _{l}$ and $%
\Gamma _{up}$ respectively, also, $\partial _{1}G$ is the vertical boundary
of $G,$ 
\begin{equation}
\Gamma _{l}=\left\{ 0<x,y<\pi ,z=B_{1}\right\} ,\text{ }\Gamma _{u}=\left\{
0<x,y<\pi ,z=B_{2}\right\} ,  \label{2.2}
\end{equation}%
\begin{equation}
\partial _{1}G=\partial G\diagdown \left( \Gamma _{l}\cup \Gamma _{u}\right)
.  \label{2.3}
\end{equation}%
Let $b\left( \mathbf{u}\right) $ be the speed of sound in $\mathbb{R}^{3}$
and let $c\left( \mathbf{u}\right) =1/b^{2}\left( \mathbf{u}\right) .$ Then $%
n\left( \mathbf{u}\right) =\sqrt{c\left( \mathbf{u}\right) }$ is the
refractive index. Let $c_{0}>0$ be a number. We assume that the following
conditions hold true:%
\begin{equation}
c\in C^{2}\left( \overline{G}\right) ,\text{ }c\in C\left( \mathbb{R}%
^{3}\right) ,  \label{2.4}
\end{equation}%
\begin{equation}
c\left( \mathbf{u}\right) \geq c_{0},\text{ }\forall \mathbf{u}\in \mathbb{R}%
^{3},  \label{2.5}
\end{equation}%
\begin{equation}
\partial _{z}c\left( \mathbf{u}\right) \geq 0,\text{ }\mathbf{u}\in G,
\label{2.6}
\end{equation}%
\begin{equation}
c\left( \mathbf{u}\right) =1\text{ for }z<B_{1}.  \label{2.7}
\end{equation}%
We remark that the monotonicity condition (\ref{2.6}) goes along well with
Geophysics, see chapter 3 in \cite{R} and \cite{Vol}. Consider the
Riemannian metric%
\begin{equation*}
dt=\sqrt{c\left( \mathbf{u}\right) }\left\vert d\mathbf{u}\right\vert ,|d%
\mathbf{u}|=\sqrt{(dx)^{2}+(dy)^{2}+(dz)^{2}}.
\end{equation*}%
For two points $\mathbf{u}$,$\mathbf{u}_{0}\in \mathbb{R}^{3}$ consider the
geodesic line $\Phi \left( \mathbf{u},\mathbf{u}_{0}\right) $ connecting
them. Then the travel time of the sound from $\mathbf{u}_{0}$ to $\mathbf{u}$
is%
\begin{equation}
t\left( \mathbf{u},\mathbf{u}_{0}\right) =\dint\limits_{\Phi \left( \mathbf{u%
},\mathbf{u}_{0}\right) }\sqrt{c\left( \mathbf{\xi }\right) }d\sigma .
\label{2.8}
\end{equation}%
For each $\mathbf{u}_{0}$ the function $t\left( \mathbf{u},\mathbf{u}%
_{0}\right) $ satisfies the eikonal equation,%
\begin{equation}
t_{z}^{2}+t_{x}^{2}+t_{y}^{2}=c\left( \mathbf{u}\right) ,  \label{2.9}
\end{equation}%
\begin{equation*}
t\left( \mathbf{u},\mathbf{u}_{0}\right) =O\left( \left\vert \mathbf{u}-%
\mathbf{u}_{0}\right\vert \right) ,\mathbf{u}\rightarrow \mathbf{u}_{0}.
\end{equation*}

We assume below that the source $\mathbf{u}_{0}$ runs along an interval $I$
of a straight line,%
\begin{equation}
\mathbf{u}_{a}\in I=\left\{ \left( x,y,z\right) :x=a\in \left[ 0,\pi \right]
,y=\pi /2,z=0\right\} .  \label{2.10}
\end{equation}%
Therefore, we use a new notation below for $t\left( \mathbf{u},\mathbf{u}%
_{0}\right) $ as $t\left( \mathbf{u},a\right) ,$ where $a\in \left[ 0,\pi %
\right] $ is the parameter in (\ref{2.10}).

We assume everywhere below the validity of the condition of the regularity
of geodesic lines, which is used in many previous publications about TTTP.

\textbf{Regularity Condition.} \emph{For every pair of points} $\left( 
\mathbf{u,u}_{a}\right) \in \overline{G}\times I$ \emph{there exists a
single geodesic line} $\Phi \left( \mathbf{u},\mathbf{u}_{a}\right) $ \emph{%
connecting them. This line intersects the boundary }$\partial G$ \emph{%
exactly twice: at a point } $s_{l}\left( \mathbf{u},\mathbf{u}_{a}\right)
\in \Gamma _{l}$\emph{\ and at another point }$\widehat{s}\left( \mathbf{u},%
\mathbf{u}_{a}\right) \in \partial G\diagdown \overline{\Gamma }_{l}$\emph{.
After reaching the point} $\widehat{s}\left( \mathbf{u},\mathbf{u}%
_{a}\right) ,$ \emph{this line leaves the domain }$G$\emph{\ and never comes
back to this domain. Also, the function} $t\left( \mathbf{u},a\right) \in
C^{2}\left( \overline{G}\times \left[ 0,\pi \right] \right) \cap C^{1}\left( 
\mathbb{R}^{3}\times \left[ 0,\pi \right] \right) .$

\textbf{Travel Time Tomography Problem (TTTP)}. \emph{Assume that conditions
(\ref{2.4})-(\ref{2.7}) hold. Determine the function} $c\left( \mathbf{u}%
\right) $ \emph{for} $\mathbf{u}\in G$ \emph{assuming that the following
function} $p\left( \mathbf{u},a\right) $ \emph{is known:}%
\begin{equation}
t\left( \mathbf{u},a\right) =p\left( \mathbf{u},a\right) ,\forall \mathbf{u}%
\in \partial G\diagdown \Gamma _{l},\forall a\in \left( 0,\pi \right) .
\label{2.11}
\end{equation}

\textbf{Remark 2.1}. \emph{It follows from (\ref{2.8}) that any geodesic
line }$\Phi \left( \mathbf{u},\mathbf{u}_{a}\right) $\emph{, which is
originated at a point }$\mathbf{u}_{a}\in I,$\emph{\ is a straight line, as
long as }$z\in \left( 0,B_{1}\right) .$\emph{\ Hence, the function }$t\left( 
\mathbf{u},a\right) $ is known for $\mathbf{u}\in \Gamma _{l},$ 
\begin{equation}
t\left( \mathbf{u},a\right) =\sqrt{\left( x-a\right) ^{2}+\left( y-\pi
/2\right) ^{2}+B_{1}^{2}}\text{ for }\mathbf{u=}\left( x,y,B_{1}\right) \in
\Gamma _{l},\text{ }a\in \left[ 0,\pi \right] .  \label{2.12}
\end{equation}

\section{Preliminaries}

\label{sec:3}

In this section we formulate/reformulate some results of \cite{KTTTP} which
we need in the follow up sections.

\subsection{A special orthonormal basis}

\label{sec:3.1}

This basis was first introduced in \cite{KJIIP}. Our numerical experience
shows that this basis works well numerically \cite{convnew,KEIT,KN}. Let $%
\beta >0$ be a number. Consider the set of functions $\left\{ \psi
_{k}\left( a\right) \right\} _{n=0}^{\infty }=\left\{ \left( a+\beta \right)
^{n}e^{a}\right\} _{n=0}^{\infty }.$ Orthonormalize these functions in the
space $L_{2}\left( 0,\pi \right) $ using the Gram-Schmidt orthonormalization
procedure. We obtain the orthonormal basis 
\begin{equation}
\left\{ \varphi _{n}\left( a\right) \right\} _{n=0}^{\infty }=\left\{
P_{n}\left( a\right) e^{a}\right\} _{n=0}^{\infty }  \label{99}
\end{equation}%
in $L_{2}\left( 0,\pi \right) .$ Here $P_{n}\left( a\right) $ is a
polynomial of the degree $n$. Note that we have done this orthonormalization
numerically in \cite{convnew,KEIT,KN} for three different inverse problems.
It works fine up to first 15 functions $\left\{ \varphi _{n}\left( a\right)
\right\} _{n=0}^{14}$ in \cite{KN}. Let $\left( \cdot ,\cdot \right) $ be
the scalar product in $L_{2}\left( 0,\pi \right) .$\ Consider the numbers $%
d_{m,n}=\left( \varphi _{m}^{\prime },\varphi _{n}\right) .$ For an integer $%
N\geq 1$ consider the $N\times N$ matrix $A_{N}=\left( d_{m,n}\right)
_{m,n=0}^{N-1}.$ Then this matrix is invertible \cite{KJIIP}.

We now present some details which were not discussed previously. Consider $%
Span\left( \varphi _{0}\left( a\right) ,...,\varphi _{N-1}\left( a\right)
\right) .$ This is an $N-$dimensional subspace $L_{2,N}\left( 0,\pi \right) $
of the space $L_{2}\left( 0,\pi \right) .$ Let $D_{N}:L_{2}\left( 0,\pi
\right) \rightarrow L_{2,N}\left( 0,\pi \right) $ be the orthogonal
projection operator. Then 
\begin{equation}
D_{N}f=\dsum\limits_{m=0}^{N-1}f_{m}\varphi _{m}\left( a\right) ,\text{ }%
f_{m}=\left( f,\varphi _{m}\right) .  \label{100}
\end{equation}%
Hence,%
\begin{equation}
\frac{d}{da}\left( D_{N}f\right) =\dsum\limits_{n=0}^{N-1}f_{m}\varphi
_{m}^{\prime }\left( a\right) .  \label{101}
\end{equation}%
It follows from (\ref{99})-(\ref{101}) that 
\begin{equation*}
\left( D_{N}f\right) \in C^{\infty }\left( 0,a\right) ,\text{ }\frac{d}{da}%
\left( D_{N}f\right) \in L_{2,N}\left( 0,\pi \right) ,\forall f\in
L_{2}\left( 0,a\right) .
\end{equation*}%
Multiply both sides of (\ref{101}) sequentially by the functions $\varphi
_{n}\left( a\right) $ and integrate with respect to $a\in \left( 0,\pi
\right) .$ Denote 
\begin{equation}
y_{n}^{N}=\left( \frac{d}{da}\left( D_{N}f\right) ,\varphi _{n}\right)
,y^{N}=\left( y_{0}^{N},...,y_{N-1}^{N}\right) ^{T},z_{N}=\left(
f_{0},...,f_{N-1}\right) ^{T}.  \label{102}
\end{equation}%
We obtain $A_{N}\left( z_{N}\right) =y^{N}.$ Hence, 
\begin{equation}
z_{N}=A_{N}^{-1}\left( y^{N}\right) .  \label{103}
\end{equation}%
It follows from (\ref{100})-(\ref{103}) that every function $\widetilde{f}%
\in L_{2,N}\left( 0,\pi \right) $ can be uniquely determined from its first
derivative without a knowledge of any initial condition $\widetilde{f}\left(
a_{0}\right) $ for any point $a_{0}\in \left( 0,\pi \right) .$ In fact, this
is the reason why this basis was originally introduced in \cite{KJIIP}.

However, even though the function $D_{N}f$ is sufficiently close to the
function $f$ in the $L_{2}\left( 0,\pi \right) -$norm for sufficiently large
values of $N$, this does not imply that functions $d\left( \left(
D_{N}f\right) \right) /da$ and $D_{N}\left( f^{\prime }\right) $ are close
to each other in the $L_{2}\left( 0,\pi \right) -$norm for $f\in C^{1}\left[
0,\pi \right] $. Furthermore, 
\begin{equation*}
\frac{d}{da}\left( D_{N}f\right) \neq D_{N}\left( f^{\prime }\right) .
\end{equation*}%
Therefore, both here and in \cite{KJIIP,convnew,KEIT,KN,KTTP}, when we
recover functions of $L_{2,N}\left( 0,\pi \right) $ from their derivatives
without knowledge of initial conditions, we work only in the approximate
sense of the orthogonal projection operator $D_{N}:L_{2}\left( 0,\pi \right)
\rightarrow L_{2,N}\left( 0,\pi \right) $.

\subsection{Estimates of functions $t_{z},t_{z}^{2}$ from the below}

\label{sec:3.2}

Let $\left( \mathbf{u},a\right) \in G\times \left( 0,\pi \right) $ be an
arbitrary pair of points and let $\Phi \left( \mathbf{u},a\right) $ be the
geodesic line connecting points $\mathbf{u}_{a}$ and $\mathbf{u}$. Let $\Phi
_{0}\left( \mathbf{u},a\right) $ be the part of this line located between
the planes $\left\{ z=0\right\} $ and $\left\{ z=B_{1}\right\} $ and let $%
\mathbf{v}_{a}=\left( x_{a},y_{a},B_{1}\right) $ be the point of the
intersection of $\Phi _{0}\left( \mathbf{u},a\right) $ with $\left\{
z=B_{1}\right\} .$ Then $0<x_{a},y_{a}<\pi $ and $\Phi _{0}\left( \mathbf{u}%
,a\right) $ is the segment of the straight line connecting points $\left(
a,\pi /2,0\right) $ and $\left( x_{a},y_{a},B_{1}\right) .$ It was shown in
the proof of Lemma 4.1 of \cite{KTTTP} that such lines can be parametrized
via the variable $z$. Hence, 
\begin{equation*}
\Phi \left( \mathbf{u},a\right) =\left\{ \left( x,y,z\right) :x=x\left(
z,a\right) ,y=y\left( z,a\right) ,z\in \left( 0,Z\left( \mathbf{u},a\right)
\right) \right\} ,
\end{equation*}%
where the number $Z=Z\left( \mathbf{u},a\right) >0$ depends on $\mathbf{u}$
and $a$ and the point $\mathbf{u}=\left( x\left( Z,a\right) ,y\left(
Z,a\right) ,Z\right) \in \partial G\diagdown \Gamma _{l}.$ It was also shown
in that proof that 
\begin{equation}
t_{z}^{2}\left( \mathbf{u},a\right) =\dint\limits_{B_{1}}^{z}c_{z}\left(
x\left( s,a\right) ,y\left( s,a\right) ,s\right) ds+t_{z}^{2}\left( \mathbf{v%
}_{a},a\right) .  \label{4.1}
\end{equation}%
By (\ref{2.12}) 
\begin{equation}
t_{z}^{2}\left( \mathbf{v}_{a},a\right) =\frac{B_{1}^{2}}{\left(
x_{a}-a\right) ^{2}+\left( y-\pi /2\right) ^{2}+B_{1}^{2}}\geq \frac{%
B_{1}^{2}}{B_{1}^{2}+5\pi ^{2}/4}.  \label{4.2}
\end{equation}

\textbf{Lemma 3.1}. \emph{Let conditions (\ref{2.4})-(\ref{2.7}) be in place.%
} \emph{Then } 
\begin{equation}
t_{z}\left( \mathbf{u},a\right) \geq \frac{B_{1}}{\sqrt{B_{1}^{2}+5\pi ^{2}/4%
}},\forall \left( \mathbf{u},a\right) \in \overline{G}\times \left[ 0,\pi %
\right] ,  \label{3.3}
\end{equation}%
\emph{\ } 
\begin{equation}
t_{z}^{2}\left( \mathbf{u},a\right) \geq \frac{B_{1}^{2}}{B_{1}^{2}+5\pi
^{2}/4},\text{ }\forall \left( \mathbf{u},a\right) \in \overline{G}\times %
\left[ 0,\pi \right] ,  \label{3.4}
\end{equation}%
\begin{equation}
t_{z}^{2}\left( \mathbf{u},a\right) -t_{z}^{2}\left( \mathbf{v}_{a},a\right)
\geq 0,\forall \left( \mathbf{u},a\right) \in \overline{G}\times \left[
0,\pi \right] .  \label{3.40}
\end{equation}

Here, (\ref{3.3}) is proven in Lemma 4.1 of \cite{KTTTP} (with slightly
different notations), (\ref{3.4}) follows from (\ref{3.3}). In addition, (%
\ref{3.4}) follows from (\ref{2.6}), (\ref{4.1}) and (\ref{4.2}). Finally, (%
\ref{3.40}) follows from (\ref{2.6}) and (\ref{4.1}).

\subsection{Eikonal equation in an integro differential form}

\label{sec:3.3}

Estimates (\ref{3.3}) and (\ref{3.4}) enable us to introduce the \emph{key
change of variables}, 
\begin{equation}
v\left( \mathbf{u},a\right) =t_{z}^{2}\left( \mathbf{u},a\right) .
\label{3.5}
\end{equation}%
By (\ref{3.3}) and (\ref{3.5}) $t_{z}\left( \mathbf{u},a\right) =\sqrt{%
v\left( \mathbf{u},a\right) }$ . Hence, using the data (\ref{2.11}), we
obtain%
\begin{equation}
t\left( x,y,z,a\right) =-\dint\limits_{z}^{B_{2}}\sqrt{v\left(
x,y,s,a\right) }ds+p\left( x,y,B_{2},a\right) .  \label{3.6}
\end{equation}%
Differentiating (\ref{3.6}), we find formulas for $t_{x}\left(
x,y,z,a\right) $ and $t_{y}\left( x,y,z,a\right) .$ Substituting these as
well as (\ref{3.5}) in (\ref{2.9}), we obtain%
\begin{equation}
v\left( x,y,z,a\right) +\left( -\dint\limits_{z}^{B_{2}}\frac{v_{x}}{2\sqrt{v%
}}\left( x,y,s,a\right) ds+p_{x}\left( x,y,B_{2},a\right) \right) ^{2}
\label{3.7}
\end{equation}%
\begin{equation*}
+\left( -\dint\limits_{z}^{B_{2}}\frac{v_{y}}{2\sqrt{v}}\left(
x,y,s,a\right) ds+p_{y}\left( x,y,B_{2},a\right) \right) ^{2}=c\left(
x,y,z\right) ,
\end{equation*}%
where $\left( x,y,z\right) \in G,a\in \left( 0,\pi \right) .$ Just as in the
first step of the Bukhgeim-Klibanov method \cite{BukhKlib,Ksurvey},
differentiate both sides of equation (\ref{3.7}) with respect to the
parameter $a$ from which the right hand side of (\ref{3.7}) does not depend.
Since $\partial _{a}c\left( \mathbf{u}\right) \equiv 0,$ then we obtain an
integro differential equation with which we work below,%
\begin{equation}
\partial _{a}v\left( x,y,z,a\right) +\frac{\partial }{\partial a}\left(
-\dint\limits_{z}^{B_{2}}\frac{v_{x}}{2\sqrt{v}}\left( x,y,s,a\right)
ds+p_{x}\left( x,y,B_{2},a\right) \right) ^{2}  \label{3.8}
\end{equation}%
\begin{equation*}
+\frac{\partial }{\partial a}\left( -\dint\limits_{z}^{B_{2}}\frac{v_{y}}{2%
\sqrt{v}}\left( x,y,s,a\right) ds+p_{y}\left( x,y,B_{2},a\right) \right)
^{2}=0,\left( x,y,z\right) \in G,a\in \left( 0,\pi \right) .
\end{equation*}

\section{Semi-finite dimensional approximation for equation (\protect\ref%
{3.8})}

\label{sec:4}

\subsection{Convenient form of the function $v\left( \mathbf{u},a\right) $}

\label{sec:4.1}

By (\ref{4.1}), (\ref{4.2}), (\ref{3.40}) and (\ref{3.5}), we seek the
solution of equation (\ref{3.8}) in the following form:%
\begin{equation}
v\left( \mathbf{u},a\right) =v\left( \mathbf{v}_{a},a\right) +\overline{v}%
\left( \mathbf{u},a\right) ,\text{ }\left( \mathbf{u},a\right) \in G\times
\left( 0,\pi \right) ,\text{ }  \label{4.3}
\end{equation}%
\begin{equation}
\overline{v}\left( \mathbf{u},a\right) \geq 0,\text{ }\left( \mathbf{u}%
,a\right) \in G\times \left( 0,\pi \right) ,  \label{4.30}
\end{equation}%
\begin{equation}
\overline{v}\left( x,y,B_{1},a\right) =0,\text{ }\left( x,y,a\right) \in
\left( 0,\pi \right) ^{3}.  \label{4.4}
\end{equation}

We need to obtain zero Dirichlet boundary condition at vertical sides of $G$%
, i.e. at $\partial _{1}G,$ see (\ref{2.3}).\ Hence, consider the data $%
p\left( \mathbf{u},a\right) $ in (\ref{2.11}). First, we apply (\ref{4.3})-(%
\ref{4.4}) to get: 
\begin{equation}
p\left( x,y,z,a\right) =p\left( x,y,B_{1},a\right) +\overline{p}\left(
x,y,z,a\right) ,\text{ }\left( x,y,z,a\right) \in \partial _{1}G\times
\left( 0,\pi \right) ,  \label{4.5}
\end{equation}%
\begin{equation}
\overline{v}\left( \mathbf{u},a\right) =\overline{p}\left( \mathbf{u}%
,a\right) ,\text{ }\left( \mathbf{u},a\right) \in \partial _{1}G\times
\left( 0,\pi \right) ,  \label{4.05}
\end{equation}%
\begin{equation}
\overline{p}\left( x,y,z,a\right) \geq 0,\text{ }\left( x,y,z,a\right) \in
\partial _{1}G\times \left( 0,\pi \right) ,  \label{4.50}
\end{equation}%
\begin{equation}
\overline{p}\left( x,y,B_{1},a\right) =0,\text{ }\left( x,y,a\right) \in
\left( 0,\pi \right) ^{3}.  \label{4.6}
\end{equation}

Based on (\ref{4.5})-(\ref{4.6}), we assume that there exists a function $%
q\left( \mathbf{u},a\right) $ such that 
\begin{equation}
q\in C^{1}\left( \overline{G}\times \left[ 0,\pi \right] \right) \text{, }%
q\left( \mathbf{u},a\right) =\overline{p}\left( \mathbf{u},a\right) ,\text{ }%
\left( \mathbf{u},a\right) \in \partial _{1}G\times \left( 0,\pi \right) ,
\label{4.7}
\end{equation}%
\begin{equation}
q\left( \mathbf{u},a\right) \geq 0,\text{ }\left( \mathbf{u},a\right) \in
\partial _{1}G\times \left( 0,\pi \right) ,\text{ }  \label{4.07}
\end{equation}%
\begin{equation}
q\left( x,y,B_{1},a\right) =0,\text{ }\left( x,y,a\right) \in \left( 0,\pi
\right) ^{3}.  \label{4.070}
\end{equation}%
We also assume that the function $q\left( \mathbf{u},a\right) $ is known.
Define the function $V\left( \mathbf{u},a\right) \in C^{1}\left( \overline{G}%
\times \left[ 0,\pi \right] \right) $ as 
\begin{equation}
V\left( \mathbf{u},a\right) =\overline{v}\left( \mathbf{u},a\right) -q\left( 
\mathbf{u},a\right) .  \label{4.70}
\end{equation}%
Then (\ref{4.4}), (\ref{4.05}), (\ref{4.6}) and (\ref{4.70}) imply that the
the function $V\left( \mathbf{u},a\right) $ satisfies the following boundary
conditions: 
\begin{equation}
V\left( \mathbf{u},a\right) =0,\text{ }\left( \mathbf{u},a\right) \in
\partial _{1}G\times \left( 0,\pi \right) ,  \label{4.8}
\end{equation}%
\begin{equation}
V\left( x,y,B_{1},a\right) =0,\text{ }\left( x,y,a\right) \in \left( 0,\pi
\right) ^{3}.  \label{4.9}
\end{equation}%
Hence, by (\ref{4.3}) and (\ref{4.70}) 
\begin{equation}
v\left( \mathbf{u},a\right) =v\left( \mathbf{v}_{a},a\right) +q\left( 
\mathbf{u},a\right) +V\left( \mathbf{u},a\right) .  \label{4.10}
\end{equation}%
In the right hand side of (\ref{4.10}), functions $v\left( \mathbf{v}%
_{a},a\right) >0$ and $q\left( \mathbf{u},a\right) $ are known, the function 
$V\left( \mathbf{u},a\right) $ is unknown and it satisfies boundary
conditions (\ref{4.8}), (\ref{4.9}). We focus below on the search of the
function $V\left( \mathbf{u},a\right) .$

\subsection{Approximate mathematical model}

\label{sec:4.2}

First, we represent the function $V\left( \mathbf{u},a\right) $ via a
truncated Fourier series with respect to the orthonormal basis $\left\{
\varphi _{k}\left( a\right) \right\} _{k=0}^{\infty }$ of section 3.1.\ More
precisely we assume that for $\left( \mathbf{u},a\right) \in \overline{G}%
\times \left[ 0,\pi \right] $ 
\begin{equation}
\partial _{x}^{k_{1}}\partial _{y}^{k_{2}}\partial _{a}^{k_{3}}V\left( 
\mathbf{u},a\right) =\dsum\limits_{n=0}^{N-1}\partial _{x}^{k_{1}}\partial
_{y}^{k_{2}}\partial _{a}^{k_{3}}\left( V_{n}\left( \mathbf{u}\right)
\varphi _{n}\left( a\right) \right) ,\text{ }\forall k_{1},k_{2},k_{3}=0,1,
\label{4.11}
\end{equation}%
where the coefficients $V_{n}\left( \mathbf{u}\right) $ are unknown. Next,
using (\ref{4.8}), we assume that for $\left( x,y,z\right) \in \left[ 0,\pi %
\right] ^{2}\times \left[ B_{1},B_{2}\right] ,n\in \left[ 0,N-1\right] $ 
\begin{equation}
V_{n}\left( \mathbf{u}\right) =\dsum\limits_{k,m=1}^{K}w_{n,km}\left(
z\right) \sin \left( kx\right) \sin \left( my\right) ,  \label{4.12}
\end{equation}%
\begin{equation}
w_{n,km}\left( z\right) \in C\left[ B_{1},B_{2}\right] ,  \label{4.13}
\end{equation}%
\begin{equation}
w_{n,km}\left( B_{1}\right) =0,  \label{4.14}
\end{equation}%
where functions $w_{n,km}\left( z\right) $ are unknown. Boundary condition (%
\ref{4.14}) is generated by (\ref{4.9}). Hence (\ref{4.10}) and (\ref{4.12})
imply that%
\begin{equation}
v\left( \mathbf{u},a\right) =v\left( \mathbf{v}_{a},a\right) +q\left( 
\mathbf{u},a\right)
+\dsum\limits_{n=0}^{N-1}\dsum\limits_{k,m=1}^{K}w_{n,km}\left( z\right)
\varphi _{n}\left( a\right) \sin \left( kx\right) \sin \left( my\right) .
\label{4.15}
\end{equation}

\textbf{Approximate Mathematical Model.} This model consists of the
assumption that the function $V\left( \mathbf{u},a\right) $ in (\ref{4.8})-(%
\ref{4.10}) can be represented via (\ref{4.11}), (\ref{4.15}) with
conditions (\ref{4.13}), (\ref{4.14}) and the substitution of so the
obtained function $v\left( \mathbf{u},a\right) $ of (\ref{4.10}) in the left
hand side of equation (\ref{3.8}) provides zero in its right hand side.
Furthermore, we assume that, in the case of noiseless data, the function $%
c^{\ast }\left( \mathbf{u}\right) $ obtained in an obvious manner from (\ref%
{4.3})-(\ref{4.15}) and (\ref{3.7}) is independent on the parameter $a$,
satisfies the first condition (\ref{2.4}) as well as conditions (\ref{2.5}),
(\ref{2.6}) in the domain $G$ and can be extended in the entire space $%
\mathbb{R}^{3}$ in such a way that it will satisfy the second condition (\ref%
{2.4}) and condition (\ref{2.7}), and, also, regularity condition holds for $%
c^{\ast }\left( \mathbf{u}\right) $.

\textbf{Remark 4.1}. \emph{The suitable values of numbers }$N$\emph{\ and }$%
K $\emph{\ should be chosen numerically, see, e.g. \cite%
{GN,Kab1,Kab2,convnew,KlibKol2,KEIT,Ktime,KN} for such choices for a variety
of inverse problems.}

Substitute (\ref{4.15}) in (\ref{3.8}). Then we obtain the following
equation for $\left( \mathbf{u},a\right) \in G\times \left( 0,\pi \right) $ 
\begin{equation*}
\dsum\limits_{k,m=1}^{K}\sin \left( kx\right) \sin \left( my\right)
\dsum\limits_{n=0}^{N-1}w_{n,km}\left( z\right) \varphi _{n}^{\prime }\left(
a\right) =
\end{equation*}%
\begin{equation}
-\frac{\partial }{\partial a}\left( -\dint\limits_{z}^{B_{2}}\frac{v_{x}}{2%
\sqrt{v}}\left( x,y,s,a\right) ds+p_{x}\left( x,y,B_{2},a\right) \right) ^{2}
\label{4.16}
\end{equation}%
\begin{equation*}
-\frac{\partial }{\partial a}\left( -\dint\limits_{z}^{B_{2}}\frac{v_{y}}{2%
\sqrt{v}}\left( x,y,s,a\right) ds+p_{y}\left( x,y,B_{2},a\right) \right)
\end{equation*}%
\begin{equation*}
-\frac{\partial }{\partial a}\left[ v\left( \mathbf{v}_{a},a\right) +q\left( 
\mathbf{u},a\right) \right] ,\text{ }v\left( \mathbf{u},a\right) \text{ is
subject to (\ref{4.12})-(\ref{4.15}).}
\end{equation*}%
Multiply both sides of equation (\ref{4.16}) by $\sin \left( kx\right) \sin
\left( my\right) ,k,m=1,...,K$ and integrate with respect to $x,y\in \left(
0,\pi \right) ^{2}.$ We obtain%
\begin{equation*}
\dsum\limits_{n=0}^{N-1}w_{n,km}\left( z\right) \varphi _{n}^{\prime }\left(
a\right) =
\end{equation*}%
\begin{equation*}
-\dint\limits_{0}^{\pi }\dint\limits_{0}^{\pi }\frac{\partial }{\partial a}%
\left( -\dint\limits_{z}^{B_{2}}\frac{v_{x}}{2\sqrt{v}}\left( x,y,s,a\right)
ds+p_{x}\left( x,y,B_{2},a\right) \right) ^{2}\sin \left( kx\right) \sin
\left( my\right) dxdy
\end{equation*}%
\begin{equation}
-\dint\limits_{0}^{\pi }\dint\limits_{0}^{\pi }\frac{\partial }{\partial a}%
\left( -\dint\limits_{z}^{B_{2}}\frac{v_{y}}{2\sqrt{v}}\left( x,y,s,a\right)
ds+p_{y}\left( x,y,B_{2},a\right) \right) ^{2}\sin \left( kx\right) \sin
\left( my\right) dxdy  \label{4.17}
\end{equation}%
\begin{equation*}
-\dint\limits_{0}^{\pi }\dint\limits_{0}^{\pi }\frac{\partial }{\partial a}%
\left[ v\left( \mathbf{u}_{a},a\right) +q\left( \mathbf{u},a\right) \right]
\sin \left( kx\right) \sin \left( my\right) dxdy,\text{ }z\in \left(
B_{1},B_{2}\right) ,
\end{equation*}%
\begin{equation*}
v\left( \mathbf{u},a\right) \text{ is subject to (\ref{4.13})-(\ref{4.15}).}
\end{equation*}%
Denote 
\begin{equation}
Q_{km}\left( z\right) =\left( w_{0,km}\left( z\right) ,...,w_{N-1,km}\left(
z\right) \right) ^{T},Q\left( z\right) =\left( Q_{km}\left( z\right) \right)
_{k,m=1}^{K}.  \label{4.18}
\end{equation}%
Multiply both sides of (\ref{4.17}) sequentially by functions $\varphi
_{0}\left( a\right) ,...,\varphi _{N-1}\left( a\right) $ and integrate then
with respect to $a\in \left( 0,\pi \right) .$ Then we obtain a system of $N$
coupled nonlinear integro differential equations with respect to the matrix $%
Q\left( z\right) .$ The left hand side of this system is $A_{N}\left(
Q_{km}\left( z\right) \right) ,$ where the invertible matrix $A_{N}$ was
introduced in section 3.1. Multiplying by $A_{N}^{-1}$ from the left and
varying $k,m=1,...,K,$ we obtain another system of coupled nonlinear integro
differential equations in the form:%
\begin{equation}
Q\left( z\right) =F\left( \dint\limits_{z}^{B_{2}}g\left( Q\left( t\right) ,%
\widehat{q}\left( t\right) \right) dt,\widehat{p}\right) +\widehat{v}_{0}+%
\widehat{q}\left( z\right) ,\text{ }z\in \left[ B_{1},B_{2}\right] ,
\label{4.19}
\end{equation}%
\begin{equation}
Q\left( B_{1}\right) =0.  \label{4.20}
\end{equation}%
Boundary condition (\ref{4.20}) follows from (\ref{4.14}). Here the vector
functions $F,g$ depend on $N,K$ and are twice continuously differentiable
with respect to their variables as long as 
\begin{equation}
q\left( \mathbf{u},a\right) +\dsum\limits_{k,m=1}^{K}\sin \left( kx\right)
\sin \left( my\right) \dsum\limits_{n=0}^{N-1}w_{n,km}\left( z\right)
\varphi _{n}\left( a\right) >-d,\text{ }\forall \left( \mathbf{u},a\right)
\in \overline{G}\times \left[ 0,\pi \right] ,  \label{4.21}
\end{equation}%
\begin{equation}
d=const\in \left( 0,\frac{B_{1}^{2}}{B_{1}^{2}+5\pi ^{2}/4}\right) .
\label{4.210}
\end{equation}%
The vector $\widehat{v}_{0}$ and the vector function $\widehat{q}\left(
z\right) \in C\left[ B_{1},B_{2}\right] $ linearly depend on functions $%
v\left( \mathbf{v}_{a},a\right) $ in (\ref{4.2}) and $q\left( \mathbf{u}%
,a\right) $ in (\ref{4.7}) respectively. The vector $\widehat{p}$ depends on
functions $p_{x}\left( x,y,B_{2},a\right) ,p_{y}\left( x,y,B_{2},a\right) .$

\textbf{Remarks 4.2:}

1\emph{.} \emph{With respect to (\ref{4.21})} \emph{recall that by (\ref%
{4.07})} $q\left( \mathbf{u},a\right) \geq 0$ \emph{in} $\overline{G}\times %
\left[ 0,\pi \right] .$ \emph{The number }$d$\emph{\ in (\ref{4.210}) is the
number of ones choice. In principle, (\ref{3.40}), (\ref{4.3}), (\ref{4.10})
and (\ref{4.12}) seem to imply that one should choose }$d=0$\emph{\ and then
replace \textquotedblleft }$>$\emph{" in (\ref{4.21}) with \textquotedblleft 
}$\geq $\emph{". However, since the function }$q\left( \mathbf{u},a\right) $ 
\emph{is generated by the data} $p\left( \mathbf{u},a\right) $, see \emph{(%
\ref{2.11}), (\ref{4.5}) and (\ref{4.7}) and since these data might contain
noise, then the noisy }$q\left( \mathbf{u},a\right) $ \emph{might be
non-positive at some points. Hence, the choice (\ref{4.210}) provides more
flexibility in terms of noise. It follows from (\ref{4.2}), (\ref{3.40}) and
(\ref{4.3}) that (\ref{4.21}) and (\ref{4.210}) guarantee that} $v\left( 
\mathbf{u},a\right) >B_{1}^{2}/\left( B_{1}^{2}+5\pi ^{2}/4\right)
-d=const.>0$ \emph{for} $\left( \mathbf{u},a\right) \in \overline{G}\times %
\left[ 0,\pi \right] .$ \emph{And we obviously need this inequality in (\ref%
{4.17}).}

2. \emph{The requirement} $v\left( \mathbf{u},a\right) \geq const.>0$ \emph{%
is a technical condition. By the numerical experience of the author in other
inverse problems, most likely, this inequality will always be satisfied in
computations if using the numerical method of this paper.}

3. \emph{Due to the boundary condition (\ref{4.20}), the problem (\ref{4.19}%
), (\ref{4.20}) cannot be solved just as a system of nonlinear coupled
Volterra integral equations.}

Let $q_{1}\left( \mathbf{u},a\right) $ and $q_{2}\left( \mathbf{u},a\right) $
be two functions $q\left( \mathbf{u},a\right) $ which generate vector
functions $\widehat{q}_{1}\left( z\right) $ and $\widehat{q}_{2}\left(
z\right) $ respectively. Then the following estimate is easy to prove%
\begin{equation}
\left\Vert \widehat{q}_{1}-\widehat{q}_{2}\right\Vert _{C\left[ B_{1},B_{2}%
\right] }\leq C_{1}\left\Vert q_{1}-q_{2}\right\Vert _{C^{1}\left( \overline{%
G}\times \left[ 0,\pi \right] \right) },  \label{4.22}
\end{equation}%
where the number $C_{1}=C_{1}\left( G,N,K\right) >0$ depends only on listed
parameters.

Here and everywhere below the norm of a vector function in a conventional
Banah space is defined as the square root of the sum of squares of norms of
its components in that space. Thus, we do not introduce below special
notations for norms of vector functions, for brevity.

Below the $NK^{2}-$dimensional vector functions $Q\left( z\right) $ have the
form (\ref{4.18}). Let $R>0$ be an arbitrary number and the number $d\in
\left( 0,B_{1}^{2}/\left( B_{1}^{2}+5\pi ^{2}/4\right) \right) $ be the one
chosen in (\ref{4.21}) (the second Remark 4.2). Denote%
\begin{equation}
H_{0}^{1}\left( B_{1},B_{2}\right) =\left\{ Q\left( z\right) \in H^{1}\left(
B_{1},B_{2}\right) :Q\left( B_{1}\right) =0\right\} ,  \label{4.23}
\end{equation}%
\begin{equation}
B\left( R,q,d\right) =\left\{ Q\left( z\right) \in H_{0}^{1}\left(
B_{1},B_{2}\right) :\left\Vert Q\right\Vert _{H^{1}\left( B_{1},B_{2}\right)
}<R,\text{ (\ref{4.21}) holds}\right\} .  \label{4.24}
\end{equation}%
Obviously, the set $B\left( R,q,d\right) $ is convex. Lemmata 4.1 and 4.2
can be proven similarly with proofs of Lemmata 5.2 and 5.3 of \cite{KTTTP}%
\emph{\ }respectively using the multidimensional analog of Taylor formula 
\cite{V}. Hence, we omit these proofs. Recall that by the embedding theorem $%
H_{0}^{1}\left( B_{1},B_{2}\right) \subset C\left[ B_{1},B_{2}\right] $
implying $B\left( R,q,d\right) \subset C\left[ B_{1},B_{2}\right] $ and 
\begin{equation}
\text{ }\left\Vert Q\right\Vert _{C\left[ B_{1},B_{2}\right] }\leq CR,\text{ 
}\forall Q\in B\left( R,q,d\right) ,  \label{4.024}
\end{equation}%
where the constant $C=C\left( N,K,B_{1},B_{2}\right) >0$ depends only on
listed parameters.

\textbf{Lemma 4.1.}\emph{\ Let }$p_{1}\left( x,y,B_{2},a\right) ,p_{2}\left(
x,y,B_{2},a\right) \in C^{1}\left( \overline{\Gamma }_{u}\times \left[ 0,\pi %
\right] \right) $\emph{\ be two functions of the data on }$\overline{\Gamma }%
_{u}$\emph{\ in (\ref{2.11}). Let} $q_{1}\left( \mathbf{u},a\right) $,$%
q_{2}\left( \mathbf{u},a\right) \in C^{1}\left( \overline{G}\times \left[
0,\pi \right] \right) $ \emph{be two functions} $q\left( \mathbf{u},a\right)
.$ \emph{For }$i=1,2$\emph{, let} 
\begin{equation}
\left\Vert p_{i}\left( x,y,B_{2},a\right) \right\Vert _{C^{1}\left( 
\overline{\Gamma }_{u}\times \left[ 0,\pi \right] \right) }\leq A,\text{ }%
\left\Vert q_{i}\left( \mathbf{u},a\right) \right\Vert _{C^{1}\left( 
\overline{G}\times \left[ 0,\pi \right] \right) }\leq A,  \label{4.240}
\end{equation}%
\emph{where }$A=const.>0.$\emph{\ Let }$\widehat{p}_{1}$\emph{\ and }$%
\widehat{p}_{2}$\emph{\ be corresponding vectors }$\widehat{p}$\emph{\
generated by }$p_{1}\left( x,y,B_{2},a\right) $\emph{\ and }$p_{2}\left(
x,y,B_{2},a\right) $\emph{\ respectively. Denote }$\widetilde{p}=\widehat{p}%
_{1}-\widehat{p}_{2}.$\emph{\ And similarly for }$\widehat{q}_{1},\widehat{q}%
_{2}$\emph{\ denote }$\widetilde{q}=\widehat{q}_{1}-\widehat{q}_{2}.$\emph{\
Then the following analog of Taylor formula is valid for any pair of vector
functions }$Q_{1}\in B\left( R,q_{1},d\right) ,Q_{2}\in B\left(
R,q_{2},d\right) $\emph{\ with }$\widetilde{Q}=Q_{2}-Q_{1}:$%
\begin{equation*}
F\left( \dint\limits_{z}^{B_{2}}g\left( Q_{2}\left( t\right) ,\widehat{q}%
_{2}\left( t\right) \right) dt,\widehat{p}_{2}\right) -F\left(
\dint\limits_{z}^{B_{2}}g\left( Q_{1}\left( t\right) ,\widehat{q}_{1}\left(
t\right) \right) dt,\widehat{p}_{1}\right)
\end{equation*}%
\begin{equation*}
=\dint\limits_{z}^{B_{2}}F_{1}\left( Q_{1}\left( t\right) ,Q_{2}\left(
t\right) ,\widehat{q}_{1}\left( t\right) ,\widehat{q}_{2}\left( t\right) ,%
\widehat{p}_{1},\widehat{p}_{2}\right) \widetilde{Q}\left( t\right) dt
\end{equation*}%
\begin{equation*}
+\dint\limits_{z}^{B_{2}}F_{2}\left( Q_{1}\left( t\right) ,Q_{2}\left(
t\right) ,\widehat{q}_{1}\left( t\right) ,\widehat{q}_{2}\left( t\right) ,%
\widehat{p}_{1},\widehat{p}_{2}\right) \widetilde{q}\left( t\right) dt
\end{equation*}%
\begin{equation*}
+F_{3}\left( Q_{1},Q_{2},\widehat{q}_{1},\widehat{q}_{2},\widehat{p}_{1},%
\widehat{p}_{2}\right) \left( z\right) \widetilde{p}
\end{equation*}%
\begin{equation*}
+\dint\limits_{z}^{B_{2}}F_{4}\left( Q_{1}\left( t\right) ,Q_{2}\left(
t\right) ,\widehat{q}_{1}\left( t\right) ,\widehat{q}_{2}\left( t\right) ,%
\widehat{p}_{1},\widehat{p}_{2},\widetilde{Q}\left( t\right) \right) dt,z\in
\left( B_{1},B_{2}\right) ,
\end{equation*}%
\emph{where vector functions }$F_{j},j=1,2,3,4$\emph{\ depend on variables
listed in them as on parameters, they are continuous with respect to these
parameters as long as for }$i=1,2$\emph{\ the vector functions }$Q_{i}\in
B\left( R,q_{i},d\right) $\emph{. The function }$F_{4}$ \emph{depends
nonlinearly on} $\widetilde{Q}.$\emph{\ Furthermore,\ the following
estimates are valid:}%
\begin{equation*}
\left\vert F_{i}\left( Q_{1}\left( t\right) ,Q_{2}\left( t\right) ,\widehat{q%
}_{1}\left( t\right) ,\widehat{q}_{2}\left( t\right) ,\widehat{p}_{1},%
\widehat{p}_{2}\right) \right\vert \leq C_{2},t\in \left[ B_{1},B_{2}\right]
,i=1,2,
\end{equation*}%
\begin{equation*}
\left\vert F_{3}\left( Q_{1},Q_{2},\widehat{q}_{1},\widehat{q}_{2},\widehat{p%
}_{1},\widehat{p}_{2}\right) \left( z\right) \right\vert \leq C_{2},\text{ }%
z\in \left[ B_{1},B_{2}\right] ,
\end{equation*}%
\begin{equation*}
\left\vert F_{4}\left( Q_{1}\left( t\right) ,Q_{2}\left( t\right) ,\widehat{q%
}_{1}\left( t\right) ,\widehat{q}_{2}\left( t\right) ,\widehat{p}_{1},%
\widehat{p}_{2},\widetilde{Q}\left( t\right) \right) \right\vert \leq
C_{2}\left\vert \widetilde{Q}\left( t\right) \right\vert ^{2},\text{ }t\in %
\left[ B_{1},B_{2}\right] .
\end{equation*}%
\emph{Here and below }$C_{2}=C_{2}\left( G,N,K,A,R,d\right) >0$\emph{\
denotes different numbers depending only on listed parameters. }

\textbf{Lemma 4.2.} \emph{Let conditions of Lemma 4.1 hold. Then}%
\begin{equation*}
F\left( \dint\limits_{z}^{B_{2}}g\left( Q_{2}\left( t\right) ,\widehat{q}%
_{2}\left( t\right) \right) dt,\widehat{p}_{2}\right) -F\left(
\dint\limits_{z}^{B_{2}}g\left( Q_{1}\left( t\right) ,\widehat{q}_{1}\left(
t\right) \right) dt,\widehat{p}_{1}\right)
\end{equation*}%
\begin{equation*}
=\dint\limits_{z}^{B_{2}}S_{1}\left( Q_{1}\left( t\right) ,Q_{2}\left(
t\right) ,\widehat{q}_{1}\left( t\right) ,\widehat{q}_{2}\left( t\right) ,%
\widehat{p}_{1},\widehat{p}_{2}\right) \widetilde{Q}\left( t\right) dt
\end{equation*}%
\begin{equation*}
+\dint\limits_{z}^{B_{2}}S_{2}\left( Q_{1}\left( t\right) ,Q_{2}\left(
t\right) ,\widehat{q}_{1}\left( t\right) ,\widehat{q}_{2}\left( t\right) ,%
\widehat{p}_{1},\widehat{p}_{2}\right) \widetilde{q}\left( t\right)
dt+S_{3}\left( Q_{1},Q_{2},\widehat{q}_{1},\widehat{q}_{2},\widehat{p}_{1},%
\widehat{p}_{2}\right) \left( z\right) \widetilde{p},
\end{equation*}%
\emph{where vector functions }$S_{k},k=1,2,3$\emph{\ have the same
properties as those of vector functions }$F_{k},k=1,2,3$\emph{\ of Lemma 4.1.%
}

\subsection{A sufficient condition ensuring (\protect\ref{4.21})}

\label{sec:4.3}

Even though (\ref{4.21}) is a technical requirement (second Remark 4.2), we
provide in this section a sufficient condition which ensures (\ref{4.21})
for $d=0$. Since our numerical method is seeking the vector function $%
Q\left( z\right) ,$ then it is convenient to formulate this condition in
terms of components $w_{n,km}\left( z\right) $ of this vector. The proof of
Lemma 4.3 is similar with the proof of Lemma 3.1 of \cite{KTTTP}. There is
an important difference, however: functions $\sin \left( kx\right) \sin
\left( my\right) $ where not used in \cite{KTTTP}.

\textbf{Lemma 4.3}. \emph{Assume that for every point }$\mathbf{u}\in 
\overline{G}$ \emph{the function} $q\left( \mathbf{u},a\right) $ \emph{in (%
\ref{4.7}) belongs to the subspace }$L_{2,N}\left( 0,\pi \right) $\emph{\
(section 3.1). Let }$V\left( \mathbf{u},a\right) $ \emph{be the function
defined in (\ref{4.70}).} \emph{Denote} $q_{n}\left( \mathbf{u}\right)
=\left( q\left( \mathbf{u},a\right) ,\varphi _{n}\left( a\right) \right) $ 
\emph{and }$V_{n}\left( \mathbf{u}\right) =\left( V\left( \mathbf{u}%
,a\right) ,\varphi _{n}\left( a\right) \right) $ \emph{where} $n=0,...,N-1.$ 
\emph{Assume that condition (\ref{4.12})} \emph{holds.} \emph{Define vector
functions }$\psi ^{N}\left( a\right) $\emph{\ and }$\varphi ^{N}\left(
a\right) $\emph{\ as }$\psi ^{N}\left( a\right) =\left( \psi _{0},...,\psi
_{N-1}\right) ^{T}\left( a\right) $\emph{\ and }$\varphi ^{N}\left( a\right)
=\left( \varphi _{0},...,\varphi _{N-1}\right) ^{T}\left( a\right) ,$\emph{\
where functions }$\psi _{n}\left( a\right) =\left( a+\beta \right) ^{n}e^{a}$%
\emph{\ were defined in section 3.1 and functions }$\varphi _{n}\left(
a\right) $\emph{\ are obtained from }$\psi ^{N}\left( a\right) $\emph{\ via
the Gram-Schmidt orthonormalization procedure, see (\ref{99}). Let }$Y_{N}$%
\emph{\ be the matrix of this procedure, }$Y_{N}\left( \psi ^{N}\right)
=\varphi ^{N}.$\emph{\ Consider the vector functions} $s\left( \mathbf{u}%
\right) =\left( q_{0}+V_{0},q_{1}+V_{1},...,q_{N-1}+V_{N-1}\right)
^{T}\left( \mathbf{u}\right) $ \emph{and} $r\left( \mathbf{u}\right)
=Y_{N}^{T}\left( s\left( \mathbf{u}\right) \right) ,$ \emph{where}

$r\left( \mathbf{u}\right) =\left( r_{0},...,r_{N-1}\right) ^{T}\left( 
\mathbf{u}\right) .$ \emph{Assume that} $r_{j}\left( \mathbf{u}\right) \geq
0 $ for all $\mathbf{u}\in \overline{G}$ \emph{and all} $j=0,...,N-1.$ \emph{%
Then} 
\begin{equation}
q\left( \mathbf{u},a\right) +V\left( \mathbf{u},a\right) \geq 0,\text{ }%
\forall \left( \mathbf{u},a\right) \in \overline{G}\times \left[ 0,\pi %
\right] .  \label{4.25}
\end{equation}

\textbf{Proof}. Let the raw number $n$ of the matrix $Y_{N}$ be $\left(
k_{n,0},k_{n,1},...,k_{n,N-1}\right) .$ Then%
\begin{equation*}
\varphi _{n}\left( a\right) =\dsum\limits_{j=0}^{N-1}k_{n,j}\psi _{j}\left(
a\right) .
\end{equation*}%
Hence,%
\begin{equation*}
q\left( \mathbf{u},a\right) +V\left( \mathbf{u},a\right)
=\dsum\limits_{n=0}^{N-1}\left[ q_{n}\left( \mathbf{u}\right) +V_{n}\left( 
\mathbf{u}\right) \right] \varphi _{n}\left( a\right)
\end{equation*}%
\begin{equation*}
=\dsum\limits_{n=0}^{N-1}\left( q_{n}\left( \mathbf{u}\right) +V_{n}\left( 
\mathbf{u}\right) \right) \dsum\limits_{j=0}^{N-1}k_{n,j}\psi _{j}\left(
a\right)
\end{equation*}%
\begin{equation*}
=\dsum\limits_{j=0}^{N-1}\left[ \dsum\limits_{n=0}^{N-1}k_{n,j}\left(
q_{n}\left( \mathbf{u}\right) +V_{n}\left( \mathbf{u}\right) \right) \right]
\psi _{j}\left( a\right) =\dsum\limits_{j=0}^{N-1}r_{j}\left( \mathbf{u}%
\right) \psi _{j}\left( a\right) .
\end{equation*}%
Since $\psi _{j}\left( a\right) >0$ for all $a\in \left[ 0,\pi \right] ,$
then (\ref{4.25}) follows. $\square $

\section{Lipschitz Stability Estimate}

\label{sec:5}

\textbf{Theorem 5.1}. \emph{Let }$R>0$\emph{\ be an arbitrary number.} For $%
i=1,2,$ \emph{let }$p^{i}\left( x,y,a\right) =p_{i}\left( x,y,B_{2},a\right)
\in C^{1}\left( \overline{\Gamma }_{u}\right) \times C\left[ 0,\pi \right] $%
\emph{\ and }$q_{i}\left( \mathbf{u},a\right) \in C^{1}\left( \overline{G}%
\right) \times C\left[ 0,\pi \right] $ \emph{be two pairs of functions
generated by the data }$p\left( \mathbf{u},a\right) $\emph{\ in (\ref{2.11})
and let estimates (\ref{4.240}) hold. Assume that for each of these pairs
there exists a solution }$Q_{i}\left( z\right) \in B\left( R,q_{i},d\right) $
\emph{of the problem (\ref{4.19}), (\ref{4.20}). Let the function }$%
c_{i}\left( \mathbf{u}\right) $ \emph{be the corresponding right hand side of%
} \emph{(\ref{3.7}). Then the following Lipschitz stability estimates are
valid: }%
\begin{equation}
\left\Vert Q_{1}-Q_{2}\right\Vert _{C\left[ B_{1},B_{2}\right] }\leq C_{2}%
\left[ \left\Vert p^{1}-p^{2}\right\Vert _{C^{1}\left( \overline{\Gamma }%
_{u}\times \left[ 0,\pi \right] \right) }+\left\Vert q_{1}-q_{2}\right\Vert
_{C^{1}\left( \overline{G}\times \left[ 0,\pi \right] \right) }\right] ,
\label{5.1}
\end{equation}%
\begin{equation}
\left\Vert c_{1}-c_{2}\right\Vert _{C\left( \overline{G}\right) }\leq C_{2}%
\left[ \left\Vert p^{1}-p^{2}\right\Vert _{C^{1}\left( \overline{\Gamma }%
_{u}\times \left[ 0,\pi \right] \right) }+\left\Vert q_{1}-q_{2}\right\Vert
_{C^{1}\left( \overline{G}\times \left[ 0,\pi \right] \right) }\right] ,
\label{5.2}
\end{equation}%
\emph{where the constant }$C_{2}=C_{2}\left( G,N,K,A,R,d\right) >0$\emph{\
was introduced in Lemma 4.1 and it depends only on listed parameters.}

\emph{\ }\textbf{Proof. }Using (\ref{4.19}) and Lemma 4.2, we obtain for $%
z\in \left[ B_{1},B_{2}\right] $ a system of linear Volterra integral
equations,%
\begin{equation}
\widetilde{Q}\left( z\right) =\dint\limits_{z}^{B_{2}}S_{1}\left(
Q_{1}\left( t\right) ,Q_{2}\left( t\right) ,\widehat{q}_{1}\left( t\right) ,%
\widehat{q}_{2}\left( t\right) ,\widehat{p}_{1},\widehat{p}_{2}\right) 
\widetilde{Q}\left( t\right) dt  \label{5.3}
\end{equation}%
\begin{equation*}
+\dint\limits_{z}^{B_{2}}S_{2}\left( Q_{1}\left( t\right) ,Q_{2}\left(
t\right) ,\widehat{q}_{1}\left( t\right) ,\widehat{q}_{2}\left( t\right) ,%
\widehat{p}_{1},\widehat{p}_{2}\right) \widetilde{q}\left( t\right)
dt+S_{3}\left( Q_{1},Q_{2},\widehat{q}_{1},\widehat{q}_{2},\widehat{p}_{1},%
\widehat{p}_{2}\right) \left( z\right) \widetilde{p}.
\end{equation*}%
Hence,%
\begin{equation*}
\left\vert \widetilde{Q}\left( z\right) \right\vert \leq C_{2}\left[
\dint\limits_{z}^{B_{2}}\widetilde{Q}\left( t\right)
dt+\dint\limits_{z}^{B_{2}}\left\vert \widetilde{q}\left( t\right)
\right\vert dt+\widetilde{p}\right] .
\end{equation*}%
Hence, (\ref{5.1}) follows from Gronwall's inequality. Finally, (\ref{5.2})
follows from (\ref{3.7}) and (\ref{5.1}). $\square $

\textbf{Remark 5.1}. \emph{Uniqueness of the reconstruction of the function} 
$c\left( \mathbf{u}\right) $ \emph{follows immediately from Lemma 5.1. In
the proof of Theorem 5.1 we have not used the boundary condition (\ref{4.20}%
) }$Q\left( B_{1}\right) =0$\emph{\ and thus came up with the system (\ref%
{5.3}) of Volterra linear integral equations. This is because Theorem 5.1 is
about the Lipschitz stability estimate rather than about a numerical method.
However, when constructing solution of equation (\ref{4.19}), we need to use
condition (\ref{4.20}). And this is done in the next section.}

\section{Numerical Method}

\label{sec:6}

In this section we construct a globally convergent numerical method which
solved problem (\ref{4.19}), (\ref{4.20}). Since the function $c\left( 
\mathbf{u}\right) $ can be straightforwardly computed from the function $%
Q\left( z\right) $ if using (\ref{3.7}), (\ref{4.3})-(\ref{4.15}), then this
method also reconstructs the target function $c\left( \mathbf{u}\right) .$

\textbf{Lemma 6.1} (Carleman estimate for the Volterra linear integral
operator) \cite{KTTTP}. \emph{Let }$\lambda >0$\emph{\ be a parameter. Then
the following Carleman estimate is valid with the Carleman Weight \ Function 
}$e^{2\lambda z}:$%
\begin{equation*}
\dint\limits_{B_{1}}^{B_{2}}\left( \dint\limits_{z}^{B_{2}}\left\vert
g\left( \tau \right) \right\vert d\tau \right) e^{2\lambda z}dz\leq \frac{1}{%
2\lambda }\dint\limits_{B_{1}}^{B_{2}}\left\vert g\left( z\right)
\right\vert e^{2\lambda z}dz,\text{ }\forall g\in L_{1}\left(
B_{1},B_{2}\right) ,\forall \lambda >0.
\end{equation*}

Let the pair of functions $\left( p\left( x,y,B_{2},a\right) ,q\left( 
\mathbf{u},a\right) \right) $ satisfies conditions (\ref{4.5}), (\ref{4.7}),
(\ref{4.070}) as well as condition (\ref{4.240}) in which $p_{1}=p_{2}=p$
and $q_{1}=q_{2}=q$. Let $d$ be an arbitrary number satisfying (\ref{4.210})
and $R>0$ be an arbitrary number. Let $B\left( R,q,d\right) \subset
H_{0}^{1}\left( B_{1},B_{2}\right) $ be the convex set defined in (\ref{4.24}%
), also see (\ref{4.23}) and (\ref{4.024}). Thus, in particular, $B\left(
R,q,d\right) $ consists of $NK^{2}-$dimensional functions $Q\left( z\right)
\in H_{0}^{1}\left( B_{1},B_{2}\right) .$

We numerically solve the problem (\ref{4.19}), (\ref{4.20}) via the
minimization on the set $B\left( R,q,d\right) $ of the following
Tikhonov-like functional with the CWF $e^{2\lambda z}$ in it$:$%
\begin{equation}
J_{\lambda ,\alpha }\left( Q\right) =e^{-2\lambda
B_{1}}\dint\limits_{B_{1}}^{B_{2}}\left[ Q-F\left(
\dint\limits_{z}^{B_{2}}g\left( Q\left( t\right) ,\widehat{q}\left( t\right)
\right) dt,\widehat{p}\right) -\widehat{v}_{0}-\widehat{q}\left( z\right) %
\right] ^{2}e^{2\lambda z}dz  \label{6.1}
\end{equation}%
\begin{equation*}
+\alpha \left\Vert Q\right\Vert _{H^{1}\left( B_{1},B_{2}\right) }^{2},\text{
}Q\in \overline{B\left( R,q,d\right) }.
\end{equation*}%
Here $\alpha \in \left( 0,1\right) $ is the regularization parameter. We use
the multiplier $e^{-2\lambda B_{1}}$ to balance two terms in the right hand
side of (\ref{6.1}). Theorem 6.1 is the main result of this paper.

\textbf{Theorem 6.1}. \emph{The functional (\ref{6.1}) has the Frech\'{e}t
derivative }$J_{\lambda ,\alpha }^{\prime }\left( Q\right) $\emph{\ at any
point of the set }$B\left( 2R,q,d\right) .$\emph{\ Furthermore, this
derivative satisfies the Lipschitz continuity condition}%
\begin{equation}
\left\Vert J_{\lambda ,\alpha }^{\prime }\left( Q_{2}\right) -J_{\lambda
,\alpha }^{\prime }\left( Q_{1}\right) \right\Vert _{H^{1}\left(
B_{1},B_{2}\right) }\leq M\left\Vert Q_{2}-Q_{1}\right\Vert _{H^{1}\left(
B_{1},B_{2}\right) },\text{ }  \label{6.2}
\end{equation}%
\emph{for all }$Q_{1},Q_{2}\in B\left( 2R,q,d\right) ,$ \emph{where the
constant }$M>0$\emph{\ is independent on }$Q_{1},Q_{2}.$ \emph{Most
importantly, there exists a sufficiently large number }%
\begin{equation*}
\overline{\lambda }=\overline{\lambda }\left( G,N,K,A,R,d\right) >1
\end{equation*}%
\emph{\ depending only on listed parameters such that the functional }$%
J_{\lambda ,\alpha }\left( Q\right) $\emph{\ is strictly convex on the set }$%
\overline{B\left( R,q,d\right) }$\emph{\ for any value of }$\lambda \geq 
\overline{\lambda }.$\emph{\ In other words, the following inequality is
valid }%
\begin{equation*}
J_{\lambda ,\alpha }\left( Q_{2}\right) -J_{\lambda ,\alpha }\left(
Q_{1}\right) -J_{\lambda ,\alpha }^{\prime }\left( Q_{1}\right) \left(
Q_{2}-Q_{1}\right)
\end{equation*}%
\begin{equation}
\geq \frac{1}{8}\left\Vert Q_{2}-Q_{1}\right\Vert _{L_{2}\left(
B_{1},B_{2}\right) }^{2}+\alpha \left\Vert Q_{2}-Q_{1}\right\Vert
_{H^{1}\left( B_{1},B_{2}\right) }^{2},\text{ }  \label{6.3}
\end{equation}%
\begin{equation*}
\forall \lambda \geq \overline{\lambda },\text{ }\forall Q_{1},Q_{2}\in 
\overline{B\left( R,q,d\right) }.
\end{equation*}

\textbf{Remark 6.1}. \emph{Although Theorem 6.1, so as other theorems of
this section, is valid only for sufficiently large values of the parameter }$%
\lambda ,$\emph{\ our past computational experience with the convexification
demonstrates that usually once can select quite reasonable values of }$%
\lambda \in \left[ 1,3\right] $\emph{\ in the numerical practice \cite%
{Bak,KlibKol1,convnew,KlibKol2,KEIT,Ktime}.}

Our next theorem is about the existence and uniqueness of the minimizer of
the functional $J_{\lambda ,\alpha }\left( Q\right) $ on the set $\overline{%
B\left( R,q,d\right) }.$

\textbf{Theorem 6.2}. \emph{Let }$\overline{\lambda }>1$\emph{\ be the
number of Theorem 6.1. Then for any }$\lambda \geq \overline{\lambda }$\emph{%
\ there exists unique minimizer }$Q_{\min }\in \overline{B\left(
R,q,d\right) }$\emph{\ of the functional }$J_{\lambda ,\alpha }\left(
Q\right) $\emph{\ on the set }$\overline{B\left( R,q,d\right) }.$\emph{\ In
addition, the following inequality is valid: }%
\begin{equation}
J_{\lambda ,\alpha }^{\prime }\left( Q_{\min }\right) \left( Q-Q_{\min
}\right) \geq 0\emph{,}\text{ }\forall Q\in \overline{B\left( R,q,d\right) }.
\label{6.30}
\end{equation}

To find the minimizer $Q_{\min }$ of Theorem 6.2, we apply the gradient
projection method. Let $T:H_{0}^{1}\left( B_{1},B_{2}\right) \rightarrow 
\overline{B\left( R,q,d\right) }$ be the orthogonal projection operator of
the space $H_{0}^{1}\left( B_{1},B_{2}\right) $ on the convex set $\overline{%
B\left( R,q,d\right) }.$ Let $\kappa \in \left( 0,1\right) $ be a number. We
choose an arbitrary point $Q^{\left( 0\right) }\in B\left( R,q,d\right) $ as
the starting point for iterations. The gradient projection method works with
the following sequence:%
\begin{equation}
Q^{\left( n\right) }=T\left( Q^{\left( n-1\right) }-\kappa J_{\lambda
,\alpha }^{\prime }\left( Q^{\left( n-1\right) }\right) \right) ,\text{ }%
n=1,2,...  \label{6.4}
\end{equation}

\textbf{Theorem 6.3}. \emph{Let }$\lambda \geq \overline{\lambda },$\emph{\
where }$\overline{\lambda }>1$\emph{\ is the number introduced in Theorem
6.1. Let }$Q_{\min }\in \overline{B\left( R,q,d\right) }$\emph{\ be the
minimizer of Theorem 6.2. Then there exists a number }$\overline{\kappa }=%
\overline{\kappa }\left( G,N,K,A,R,d,\lambda \right) \in \left( 0,1\right) $%
\emph{\ depending only on listed parameters such that for every number }$%
\kappa \in \left( 0,\overline{\kappa }\right) $\emph{\ there exists a number 
}$\theta =\theta \left( \kappa \right) \in \left( 0,1\right) $\emph{\ such
that} 
\begin{equation*}
\left\Vert Q^{\left( n\right) }-Q_{\min }\right\Vert _{H^{1}\left(
B_{1},B_{2}\right) }\leq \theta ^{n}\left\Vert Q^{\left( 0\right) }-Q_{\min
}\right\Vert _{H^{1}\left( B_{1},B_{2}\right) }.
\end{equation*}

To prove the convergence of the sequence (\ref{6.4}) to the exact solution,
we assume first that there exists exact solution $Q^{\ast }\in B\left(
R,q^{\ast },d\right) $ of the problem (\ref{4.19}), (\ref{4.20}) with
idealized noiseless data $\left( p^{\ast }\left( x,y,B_{2},a\right) ,q^{\ast
}\left( \mathbf{u},a\right) \right) .$ Such an assumption is a common place
in the theory of ill-posed problems, see, e.g. \cite{BK,T}. We also assume
that functions $p^{\ast },q^{\ast }$ satisfy conditions (\ref{4.5})-(\ref%
{4.070}) as well as condition (\ref{4.240}) in which $p_{1}=p_{2}=p^{\ast }$
and $q_{1}=q_{2}=q^{\ast }.$ Next, let $\delta \in \left( 0,d/2\right) $ be
the level of noise in the data. Hence, we assume that 
\begin{equation}
\left\Vert \left( p-p^{\ast }\right) \left( x,y,B_{2},a\right) \right\Vert
_{C^{1}\left( \overline{\Gamma }_{u}\times \left[ 0,\pi \right] \right) },%
\text{ }\left\Vert q-q^{\ast }\right\Vert _{C^{1}\left( \overline{G}\times %
\left[ 0,\pi \right] \right) }\leq \delta .  \label{6.5}
\end{equation}

Let $c^{\ast }\left( \mathbf{u}\right) $ be the exact coefficient $c\left( 
\mathbf{u}\right) $ which corresponds to the above exact data within the
framework of our approximate mathematical model of section 4.2. Combining (%
\ref{4.3})-(\ref{4.15}) with (\ref{3.7}) in an obvious manner, one can
construct approximations $\overline{c}^{\left( n\right) }\left( \mathbf{u}%
,a\right) $ for the target coefficient $c^{\ast }\left( \mathbf{u}\right) .$
Since functions $\overline{c}^{\left( n\right) }\left( \mathbf{u},a\right) $
might depend on the parameter $a\in \left[ 0,\pi \right] ,$ then we average
these approximations with respect to $a$, thus obtaining functions $%
c^{\left( n\right) }\left( \mathbf{u}\right) ,$%
\begin{equation*}
c^{\left( n\right) }\left( \mathbf{u}\right) =\frac{1}{\pi }%
\dint\limits_{0}^{\pi }\overline{c}^{\left( n\right) }\left( \mathbf{u}%
,a\right) da.
\end{equation*}

\textbf{Theorem 6.4}. \emph{Assume that conditions of Theorem 6.3 hold. Fix
an arbitrary number }$\widetilde{\lambda }\geq \overline{\lambda }$\emph{\
and set }$\lambda =\widetilde{\lambda }$ \emph{in the functional }$%
J_{\lambda ,\alpha }\left( Q\right) .$\emph{\ Also, let the regularization
parameter }$\alpha =\alpha \left( \delta \right) =\delta ^{2}.$\emph{\ Then
the following convergence estimates are valid:}%
\begin{equation}
\left\Vert Q^{\ast }-Q_{\min }\right\Vert _{L_{2}\left( B_{1},B_{2}\right)
}\leq C_{2}\delta ,  \label{6.6}
\end{equation}%
\begin{equation}
\left\Vert Q^{\ast }-Q^{\left( n\right) }\right\Vert _{L_{2}\left(
B_{1},B_{2}\right) }\leq C_{2}\delta +\theta ^{n}\left\Vert Q^{\left(
0\right) }-Q_{\min }\right\Vert _{H^{1}\left( B_{1},B_{2}\right) },
\label{6.7}
\end{equation}%
\begin{equation}
\left\Vert c^{\ast }-c^{\left( n\right) }\right\Vert _{L_{2}\left(
B_{1},B_{2}\right) }\leq C_{2}\delta +\theta ^{n}\left\Vert Q^{\left(
0\right) }-Q_{\min }\right\Vert _{H^{1}\left( B_{1},B_{2}\right) }.
\label{6.8}
\end{equation}

\textbf{Remark 6.2}.\emph{\ Since }$Q^{\left( 0\right) }$\emph{\ is an
arbitrary point of the set }$B\left( R,q,d\right) $\emph{\ and since the
size of this set }$R>0$\emph{\ is an arbitrary fixed number, then estimates (%
\ref{6.6})-(\ref{6.8}) mean global convergence, as defined in the first
paragraph of Introduction.}

As soon as Theorem 6.1 is established, lemma 2.1 of \cite{Bak} implies
Theorem 6.2. Next, if Theorems 6.1 and 6.2 are true, then theorem 2.1 \cite%
{Bak} leads to Theorem 6.3. Therefore, we are left to prove only two
theorems: 6.1 and 6.4.

\section{Proving Theorems 6.1 and 6.4}

\label{sec:7}

\subsection{Theorem 6.1}

\label{sec:7.1}

The goal of this section is to prove Theorem 6.1. Let $Q_{1},Q_{2}\in
B\left( R,q,d\right) $ be two arbitrary points of this set. Denote $%
\widetilde{Q}=Q_{2}-Q_{1}\in H_{0}^{1}\left( B_{1},B_{2}\right) .$ Then $%
Q_{2}=Q_{1}+\widetilde{Q}$ and $\left\Vert \widetilde{Q}\right\Vert
_{H^{1}\left( B_{1},B_{2}\right) }<2R.$ Hence, by (\ref{4.024})%
\begin{equation}
\left\Vert Q_{1}\right\Vert _{C\left[ B_{1},B_{2}\right] },\text{ }%
\left\Vert Q_{2}\right\Vert _{C\left[ B_{1},B_{2}\right] },\text{ }%
\left\Vert \widetilde{Q}\right\Vert _{C\left[ B_{1},B_{2}\right] }<CR.
\label{7.1}
\end{equation}%
Setting in Lemma 4.1 $\widehat{p}_{1}=\widehat{p}_{2}=\widehat{p}$ and $%
\widehat{q}_{1}=\widehat{q}_{2}=\widehat{q}$, we obtain for $z\in \left(
B_{1},B_{2}\right) $%
\begin{equation*}
Q_{2}\left( z\right) -F\left( \dint\limits_{z}^{B_{2}}g\left( Q_{2}\left(
t\right) ,\widehat{q}\left( t\right) \right) dt,\widehat{p}\right) -\widehat{%
v}_{0}-\widehat{q}\left( z\right)
\end{equation*}%
\begin{equation*}
=Q_{1}\left( z\right) -F\left( \dint\limits_{z}^{B_{2}}g\left( Q_{1}\left(
t\right) ,\widehat{q}\left( t\right) \right) dt,\widehat{p}\right) -\widehat{%
v}_{0}-\widehat{q}\left( z\right)
\end{equation*}%
\begin{equation*}
+\widetilde{Q}\left( z\right) +\dint\limits_{z}^{B_{2}}D_{1}\left(
Q_{1}\left( t\right) ,Q_{2}\left( t\right) ,\widehat{q}\left( t\right) ,%
\widehat{p}\right) \widetilde{Q}\left( t\right) dt
\end{equation*}%
\begin{equation*}
+\dint\limits_{z}^{B_{2}}D_{2}\left( Q_{1}\left( t\right) ,Q_{2}\left(
t\right) ,\widehat{q}\left( t\right) ,\widehat{p},\widetilde{Q}\left(
t\right) \right) dt,
\end{equation*}%
where $D_{1}$ and $D_{2}$ are continuos functions of their variables, the
function $D_{2}$ depends nonlinearly on $\widetilde{Q}$ and the following
two estimates are valid:%
\begin{equation}
\left\vert D_{1}\left( Q_{1}\left( t\right) ,Q_{2}\left( t\right) ,\widehat{q%
}\left( t\right) ,\widehat{p}\right) \right\vert \leq C_{2},\text{ \ }t\in %
\left[ B_{1},B_{2}\right] ,  \label{7.2}
\end{equation}%
\begin{equation}
\left\vert D_{2}\left( Q_{1}\left( t\right) ,Q_{2}\left( t\right) ,\widehat{q%
}\left( t\right) ,\widehat{p},\widetilde{Q}\left( t\right) \right)
\right\vert \leq C_{2}\left\vert \widetilde{Q}\left( t\right) \right\vert
^{2},\text{ \ }t\in \left[ B_{1},B_{2}\right] .  \label{7.3}
\end{equation}%
Denote 
\begin{equation}
X\left( z\right) =Q_{1}\left( z\right) -F\left(
\dint\limits_{z}^{B_{2}}g\left( Q_{1}\left( t\right) ,\widehat{q}\left(
t\right) \right) dt,\widehat{p}\right) -\widehat{v}_{0}-\widehat{q}\left(
z\right) .  \label{7.4}
\end{equation}%
Hence,%
\begin{equation*}
\left[ Q_{2}\left( z\right) -F\left( \dint\limits_{z}^{B_{2}}g\left(
Q_{2}\left( t\right) ,\widehat{q}\left( t\right) \right) dt,\widehat{p}%
\right) -\widehat{v}_{0}-\widehat{q}\left( z\right) \right] ^{2}-X^{2}\left(
z\right)
\end{equation*}%
\begin{equation*}
=2X\left( z\right) \left[ \widetilde{Q}\left( z\right)
+\dint\limits_{z}^{B_{2}}D_{1}\left( Q_{1}\left( t\right) ,Q_{2}\left(
t\right) ,\widehat{q}\left( t\right) ,\widehat{p}\right) \widetilde{Q}\left(
t\right) dt\right]
\end{equation*}%
\begin{equation*}
+2X\left( z\right) \dint\limits_{z}^{B_{2}}D_{2}\left( Q_{1}\left( t\right)
,Q_{2}\left( t\right) ,\widehat{q}\left( t\right) ,\widehat{p},\widetilde{Q}%
\left( t\right) \right) dt
\end{equation*}%
\begin{equation}
+\left[ \widetilde{Q}\left( z\right) +\dint\limits_{z}^{B_{2}}D_{1}\left(
Q_{1}\left( t\right) ,Q_{2}\left( t\right) ,\widehat{q}\left( t\right) ,%
\widehat{p}\right) \widetilde{Q}\left( t\right) dt\right] ^{2}  \label{7.5}
\end{equation}%
\begin{equation*}
+2\left[ \widetilde{Q}\left( z\right) +\dint\limits_{z}^{B_{2}}D_{1}\left(
Q_{1}\left( t\right) ,Q_{2}\left( t\right) ,\widehat{q}\left( t\right) ,%
\widehat{p}\right) \widetilde{Q}\left( t\right) dt\right]
\end{equation*}%
\begin{equation*}
\times \left[ \dint\limits_{z}^{B_{2}}D_{2}\left( Q_{1}\left( t\right)
,Q_{2}\left( t\right) ,\widehat{q}\left( t\right) ,\widehat{p},\widetilde{Q}%
\left( t\right) \right) dt\right]
\end{equation*}%
\begin{equation*}
+\left[ \dint\limits_{z}^{B_{2}}D_{2}\left( Q_{1}\left( t\right)
,Q_{2}\left( t\right) ,\widehat{q}\left( t\right) ,\widehat{p},\widetilde{Q}%
\left( t\right) \right) dt\right] ^{2}.
\end{equation*}%
First, we estimate from the below all lines of (\ref{7.5}), except of the
first two. Using (\ref{7.1}), (\ref{7.3}) and (\ref{7.4}), we obtain for the
third line:%
\begin{equation}
2X\left( z\right) \dint\limits_{z}^{B_{2}}D_{2}\left( Q_{1}\left( t\right)
,Q_{2}\left( t\right) ,\widehat{q}\left( t\right) ,\widehat{p},\widetilde{Q}%
\left( t\right) \right) dt\geq -C_{2}\dint\limits_{z}^{B_{2}}\left\vert 
\widetilde{Q}\left( t\right) \right\vert ^{2}dt.  \label{7.6}
\end{equation}%
Next, using Cauchy-Schwarz inequality, we obtain for the fourth line:%
\begin{equation}
\left[ \widetilde{Q}\left( z\right) +\dint\limits_{z}^{B_{2}}D_{1}\left(
Q_{1}\left( t\right) ,Q_{2}\left( t\right) ,\widehat{q}\left( t\right) ,%
\widehat{p}\right) \widetilde{Q}\left( t\right) dt\right] ^{2}\geq \frac{1}{2%
}\left\vert \widetilde{Q}\left( z\right) \right\vert
^{2}-C_{2}\dint\limits_{z}^{B_{2}}\left\vert \widetilde{Q}\left( t\right)
\right\vert ^{2}dt.  \label{7.7}
\end{equation}%
Similarly for the fifth, sixth and seventh lines,%
\begin{equation}
+2\left[ \widetilde{Q}\left( z\right) +\dint\limits_{z}^{B_{2}}D_{1}\left(
Q_{1}\left( t\right) ,Q_{2}\left( t\right) ,\widehat{q}\left( t\right) ,%
\widehat{p}\right) \widetilde{Q}\left( t\right) dt\right]  \label{7.8}
\end{equation}%
\begin{equation*}
\times \left[ \dint\limits_{z}^{B_{2}}D_{2}\left( Q_{1}\left( t\right)
,Q_{2}\left( t\right) ,\widehat{q}\left( t\right) ,\widehat{p},\widetilde{Q}%
\left( t\right) \right) dt\right] \geq -\frac{1}{4}\left\vert \widetilde{Q}%
\left( z\right) \right\vert ^{2}-C_{2}\dint\limits_{z}^{B_{2}}\left\vert 
\widetilde{Q}\left( t\right) \right\vert ^{2}dt,
\end{equation*}%
\begin{equation}
\left[ \dint\limits_{z}^{B_{2}}D_{2}\left( Q_{1}\left( t\right) ,Q_{2}\left(
t\right) ,\widehat{q}\left( t\right) ,\widehat{p},\widetilde{Q}\left(
t\right) \right) dt\right] ^{2}\geq 0.  \label{7.9}
\end{equation}%
Denote $Z$ the sum of all lines of (\ref{7.5}), except of the first two.
Hence, (\ref{7.6})-(\ref{7.9}) imply%
\begin{equation}
Z\geq \frac{1}{4}\left\vert \widetilde{Q}\left( z\right) \right\vert
^{2}-C_{2}\dint\limits_{z}^{B_{2}}\left\vert \widetilde{Q}\left( t\right)
\right\vert ^{2}dt.  \label{7.10}
\end{equation}

The second line of (\ref{7.5}) is linear with respect to $\widetilde{Q}%
\left( z\right) .$ Hence, (\ref{6.1}) and (\ref{7.5}) lead to%
\begin{equation}
J_{\lambda ,\alpha }\left( Q_{2}\right) -J_{\lambda ,\alpha }\left(
Q_{1}\right) =I_{lin}\left( \widetilde{Q}\right) +I_{nonlin}\left( 
\widetilde{Q}\right) ,  \label{7.11}
\end{equation}%
where $I_{lin}\left( \widetilde{Q}\right) $ depends linearly and $%
I_{nonlin}\left( \widetilde{Q}\right) $ depends nonlinearly on $\widetilde{Q}%
.$ \ More precisely, the linear part is:%
\begin{equation*}
I_{lin}\left( \widetilde{Q}\right) =
\end{equation*}%
\begin{equation}
=2e^{-2\lambda B_{1}}\dint\limits_{B_{1}}^{B_{2}}X\left( z\right) \left[ 
\widetilde{Q}\left( z\right) +\dint\limits_{z}^{B_{2}}D_{1}\left(
Q_{1}\left( t\right) ,Q_{2}\left( t\right) ,\widehat{q}\left( t\right) ,%
\widehat{p}\right) \widetilde{Q}\left( t\right) dt\right] e^{2\lambda z}dz
\label{7.12}
\end{equation}%
\begin{equation*}
+2\alpha \left[ Q_{1},\widetilde{Q}\right] ,
\end{equation*}%
where $\left[ \cdot ,\cdot \right] $ denotes the scalar product in the space 
$H^{1}\left( B_{1},B_{2}\right) $ of $NK^{2}-$dimensional functions. As to
the nonlinear part, it satisfies estimates (\ref{7.5}), (\ref{7.14}), where (%
\ref{7.13}) is obvious and (\ref{7.14}) follows from (\ref{7.10}): 
\begin{equation}
\left\vert I_{nonlin}\left( \widetilde{Q}\right) \right\vert \leq
C_{2}\left\Vert \widetilde{Q}\right\Vert _{C\left[ B_{1},B_{2}\right]
}^{2}\exp \left[ 2\lambda \left( B_{2}-B_{1}\right) \right] ,  \label{7.13}
\end{equation}%
\begin{equation}
I_{nonlin}\left( \widetilde{Q}\right) \geq \frac{1}{4}e^{-2\lambda
B_{1}}\dint\limits_{B_{1}}^{B_{2}}\left\vert \widetilde{Q}\left( z\right)
\right\vert ^{2}e^{2\lambda z}dz-C_{2}e^{-2\lambda
B_{1}}\dint\limits_{B_{1}}^{B_{2}}\left( \dint\limits_{z}^{B_{2}}\left\vert 
\widetilde{Q}\left( t\right) \right\vert ^{2}dt\right) e^{2\lambda z}dz
\label{7.14}
\end{equation}%
\begin{equation*}
+\alpha \left\Vert \widetilde{Q}\right\Vert _{H^{1}\left( B_{1},B_{2}\right)
}^{2}.
\end{equation*}

It follows from (\ref{7.2}), (\ref{7.4}) and (\ref{7.12}) that 
\begin{equation*}
I_{lin}\left( P\right) \leq C_{2}\exp \left[ 2\lambda \left(
B_{2}-B_{1}\right) \right] \left\Vert P\right\Vert _{H_{0}^{1}\left(
B_{1},B_{2}\right) },\text{ }\forall P\in H_{0}^{1}\left( B_{1},B_{2}\right)
.
\end{equation*}%
Hence, by Riesz theorem there exists unique element $Y\in H_{0}^{1}\left(
B_{1},B_{2}\right) $ such that 
\begin{equation}
I_{lin}\left( P\right) =\left[ Y,P\right] ,\text{ }\forall P\in
H_{0}^{1}\left( B_{1},B_{2}\right) .  \label{7.15}
\end{equation}%
In addition, it follows from (\ref{7.11}), (\ref{7.13}) and (\ref{7.15})
that 
\begin{equation*}
\lim_{\left\Vert \widetilde{Q}\right\Vert _{H_{0}^{1}\left(
B_{1},B_{2}\right) }\rightarrow 0}\left\{ \frac{1}{\left\Vert \widetilde{Q}%
\right\Vert _{H_{0}^{1}\left( B_{1},B_{2}\right) }}\left[ J_{\lambda ,\alpha
}\left( Q_{2}\right) -J_{\lambda ,\alpha }\left( Q_{1}\right) -\left[ Y,%
\widetilde{Q}\right] \right] \right\} =0.
\end{equation*}%
Hence, $I_{lin}\left( \widetilde{Q}\right) =\left[ Y,\widetilde{Q}\right] $
is the Frech\'{e}t derivative of the functional $J_{\lambda ,\alpha }\left(
Q\right) $ at the point $Q_{1},$ i.e. 
\begin{equation}
I_{lin}\left( \widetilde{Q}\right) =J_{\lambda ,\alpha }^{\prime }\left(
Q_{1}\right) \left( \widetilde{Q}\right) .  \label{7.16}
\end{equation}%
We omit the proof of the Lipschitz continuity property (\ref{6.2}) of the
Frech\'{e}t derivative $J_{\lambda ,\alpha }^{\prime }\left( Q\right) $
since this proof completely similar with the proof of Theorem 3.1 of \cite%
{Bak}.

It follows from (\ref{7.11}), (\ref{7.14}) and (\ref{7.16}) that 
\begin{equation*}
J_{\lambda ,\alpha }\left( Q_{1}+\widetilde{Q}\right) -J_{\lambda ,\alpha
}\left( Q_{1}\right) -J_{\lambda ,\alpha }^{\prime }\left( Q_{1}\right)
\left( \widetilde{Q}\right)
\end{equation*}%
\begin{equation}
\geq \frac{1}{4}e^{-2\lambda B_{1}}\dint\limits_{B_{1}}^{B_{2}}\left\vert 
\widetilde{Q}\left( z\right) \right\vert ^{2}e^{2\lambda
z}dz-C_{2}e^{-2\lambda B_{1}}\dint\limits_{B_{1}}^{B_{2}}\left(
\dint\limits_{z}^{B_{2}}\left\vert \widetilde{Q}\left( t\right) \right\vert
^{2}dt\right) e^{2\lambda z}dz  \label{7.17}
\end{equation}%
\begin{equation*}
+\alpha \left\Vert \widetilde{Q}\right\Vert _{H^{1}\left( B_{1},B_{2}\right)
}^{2}.
\end{equation*}%
By Lemma 6.1%
\begin{equation*}
C_{2}e^{-2\lambda B_{1}}\dint\limits_{B_{1}}^{B_{2}}\left(
\dint\limits_{z}^{B_{2}}\left\vert \widetilde{Q}\left( t\right) \right\vert
^{2}dt\right) e^{2\lambda z}dz\leq \frac{1}{2\lambda }C_{2}e^{-2\lambda
B_{1}}\dint\limits_{B_{1}}^{B_{2}}\left\vert \widetilde{Q}\left( z\right)
\right\vert ^{2}e^{2\lambda z}dz.
\end{equation*}%
Choose $\overline{\lambda }>1$ so large that $C_{2}/\left( 2\lambda \right)
<1/8$ for all $\lambda \geq \overline{\lambda }.$ Hence, (\ref{7.17})
implies that 
\begin{equation*}
J_{\lambda ,\alpha }\left( Q_{1}+\widetilde{Q}\right) -J_{\lambda ,\alpha
}\left( Q_{1}\right) -J_{\lambda ,\alpha }^{\prime }\left( Q_{1}\right)
\left( \widetilde{Q}\right)
\end{equation*}%
\begin{equation*}
\geq \frac{1}{8}\left\Vert \widetilde{Q}\right\Vert _{L_{2}\left(
B_{1},B_{2}\right) }^{2}+\alpha \left\Vert \widetilde{Q}\right\Vert
_{H^{1}\left( B_{1},B_{2}\right) }^{2}.\text{ }\square
\end{equation*}

\subsection{Theorem 6.4}

\label{sec:7.2}

In this section we prove Theorem 6.4. Setting in Lemma 4.2 
\begin{equation}
p_{1}\left( x,y,B_{2},a\right) =p^{\ast }\left( x,y,B_{2},a\right)
,p_{2}\left( x,y,B_{2},a\right) =p\left( x,y,B_{2},a\right) ,  \label{7.18}
\end{equation}%
\begin{equation}
q_{1}\left( \mathbf{u},a\right) =q^{\ast }\left( \mathbf{u},a\right)
,q_{2}\left( \mathbf{u},a\right) =q\left( \mathbf{u},a\right) ,Q_{1}\left(
z\right) =Q_{2}\left( z\right) =Q^{\ast }\left( z\right) ,  \label{7.19}
\end{equation}%
we obtain 
\begin{equation*}
Q^{\ast }\left( z\right) -F\left( \dint\limits_{z}^{B_{2}}g\left( Q^{\ast
}\left( t\right) ,\widehat{q}\left( t\right) \right) dt,\widehat{p}\right) -%
\widehat{v}_{0}-\widehat{q}\left( z\right)
\end{equation*}%
\begin{equation*}
=Q^{\ast }\left( z\right) -F\left( \dint\limits_{z}^{B_{2}}g\left( Q^{\ast
}\left( t\right) ,\widehat{q}^{\ast }\left( t\right) \right) dt,\widehat{p}%
^{\ast }\right) -\widehat{v}_{0}-\widehat{q}^{\ast }\left( z\right)
\end{equation*}%
\begin{equation}
-\dint\limits_{z}^{B_{2}}S_{2}\left( Q^{\ast }\left( t\right) ,Q^{\ast
}\left( t\right) ,\widehat{q}^{\ast }\left( t\right) ,\widehat{q}\left(
t\right) ,\widehat{p}^{\ast },\widehat{p}\right) \widetilde{q}\left(
t\right) dt  \label{7.20}
\end{equation}%
\begin{equation*}
-S_{3}\left( Q^{\ast }\left( t\right) ,Q^{\ast }\left( t\right) ,\widehat{q}%
^{\ast }\left( t\right) ,\widehat{q}\left( t\right) ,\widehat{p}^{\ast },%
\widehat{p}\right) \left( z\right) \widetilde{p}-\widetilde{q}\left(
z\right) .
\end{equation*}%
Since 
\begin{equation*}
Q^{\ast }\left( z\right) -F\left( \dint\limits_{z}^{B_{2}}g\left( Q^{\ast
}\left( t\right) ,\widehat{q}^{\ast }\left( t\right) \right) dt,\widehat{p}%
^{\ast }\right) -\widehat{v}_{0}-\widehat{q}^{\ast }\left( z\right) =0,
\end{equation*}%
then (\ref{6.5}) and (\ref{7.18})-(\ref{7.20}) imply that 
\begin{equation}
\left( Q^{\ast }\left( z\right) -F\left( \dint\limits_{z}^{B_{2}}g\left(
Q^{\ast }\left( t\right) ,\widehat{q}\left( t\right) \right) dt,\widehat{p}%
\right) -\widehat{v}_{0}-\widehat{q}\left( z\right) \right) ^{2}\leq
C_{2}\delta ^{2}.  \label{7.21}
\end{equation}

To indicate the dependence of the functional (\ref{6.1}) on vector functions 
$p\left( x,y,B_{2},a\right) $ and $q\left( \mathbf{u},a\right) ,$ we
temporary denote this functional as $J_{\lambda ,\alpha }\left( Q,\widehat{p}%
,\widehat{q}\right) .$ Recall that $\alpha =\delta ^{2}.$ Since the number $%
\lambda =\widetilde{\lambda }\geq \overline{\lambda }$ is not varied, then
we can estimate $\exp \left[ 2\lambda \left( B_{2}-B_{1}\right) \right] \leq
C_{2}.$ Hence, using (\ref{6.1}) and (\ref{7.21}), we obtain with a
different constant $C_{2}>0:$ 
\begin{equation}
J_{\lambda ,\alpha }\left( Q^{\ast },\widehat{p},\widehat{q}\right) \leq
C_{2}\delta ^{2}.  \label{7.22}
\end{equation}%
Next, by (\ref{6.3})%
\begin{equation*}
J_{\lambda ,\alpha }\left( Q^{\ast },\widehat{p},\widehat{q}\right)
-J_{\lambda ,\alpha }\left( Q_{\min },\widehat{p},\widehat{q}\right)
-J_{\lambda ,\alpha }^{\prime }\left( Q_{\min },\widehat{p},\widehat{q}%
\right) \left( Q^{\ast }-Q_{\min }\right) 
\end{equation*}%
\begin{equation*}
\geq \frac{1}{8}\left\Vert Q^{\ast }-Q_{\min }\right\Vert _{L_{2}\left(
B_{1},B_{2}\right) }^{2}.
\end{equation*}%
Hence, using (\ref{7.22}), we obtain%
\begin{equation}
\left\Vert Q^{\ast }-Q_{\min }\right\Vert _{L_{2}\left( B_{1},B_{2}\right)
}^{2}\leq C_{2}\delta ^{2}-J_{\lambda ,\alpha }^{\prime }\left( Q_{\min },%
\widehat{p},\widehat{q}\right) \left( Q^{\ast }-Q_{\min }\right) .
\label{7.23}
\end{equation}%
By (\ref{6.30}) 
\begin{equation}
-J_{\lambda ,\alpha }^{\prime }\left( Q_{\min },\widehat{p},\widehat{q}%
\right) \left( Q^{\ast }-Q_{\min }\right) \leq 0.  \label{7.24}
\end{equation}%
Using (\ref{7.23}) and (\ref{7.24}) we obtain (\ref{6.6}).

Next, the triangle inequality and Theorem 6.3 give:%
\begin{equation*}
\left\Vert Q^{\ast }-Q^{\left( n\right) }\right\Vert _{L_{2}\left(
B_{1},B_{2}\right) }\leq \left\Vert Q^{\ast }-Q_{\min }\right\Vert
_{L_{2}\left( B_{1},B_{2}\right) }+\left\Vert Q^{\left( n\right) }-Q_{\min
}\right\Vert _{L_{2}\left( B_{1},B_{2}\right) }
\end{equation*}%
\begin{equation*}
\leq \left\Vert Q^{\ast }-Q_{\min }\right\Vert _{L_{2}\left(
B_{1},B_{2}\right) }+\left\Vert Q^{\left( n\right) }-Q_{\min }\right\Vert
_{H^{1}\left( B_{1},B_{2}\right) }
\end{equation*}%
\begin{equation*}
\leq \left\Vert Q^{\ast }-Q_{\min }\right\Vert _{L_{2}\left(
B_{1},B_{2}\right) }+\theta ^{n}\left\Vert Q^{\left( 0\right) }-Q_{\min
}\right\Vert _{H^{1}\left( B_{1},B_{2}\right) }.
\end{equation*}%
Thus, 
\begin{equation}
\left\Vert Q^{\ast }-Q^{\left( n\right) }\right\Vert _{L_{2}\left(
B_{1},B_{2}\right) }\leq \left\Vert Q^{\ast }-Q_{\min }\right\Vert
_{L_{2}\left( B_{1},B_{2}\right) }+\theta ^{n}\left\Vert Q^{\left( 0\right)
}-Q_{\min }\right\Vert _{H^{1}\left( B_{1},B_{2}\right) }.  \label{7.25}
\end{equation}%
Combining (\ref{7.25}) with (\ref{6.6}), we obtain (\ref{6.7}). Finally,
estimate (\ref{6.8}) follows (\ref{6.7}) and the construction of functions $%
c^{\left( n\right) }\left( \mathbf{u}\right) $ described in section 6. $%
\square $

\begin{center}
\textbf{Acknowledgment}
\end{center}

This work was supported by US Army Research Laboratory and US Army Research
Office grant W911NF-19-1-0044.

\end{document}